\documentclass[reqno,11pt]{amsart}

\usepackage{amsmath, amsthm, amssymb, amsfonts, amsthm}
\usepackage{amsxtra, amscd, geometry, graphicx,epic,eepic}
\usepackage[all,cmtip]{xy}
\usepackage{mathrsfs}
\usepackage{hyperref}

\theoremstyle{definition}

\newtheorem*{rmk*}{Remark}
\newtheorem*{thm1.1*}{Theorem 1.1}
\newtheorem*{thm1.2*}{Theorem 1.2}
\newtheorem*{cor1.3*}{Corollary 1.3}
\newtheorem*{lem2.1*}{Lemma 2.1}
\newtheorem*{cor2.2*}{Corollary 2.2}
\newtheorem*{prop2.3*}{Proposition 2.3}
\newtheorem*{rem2.4*}{Remark 2.4}
\newtheorem*{prop2.5*}{Proposition 2.5}
\newtheorem*{lem2.6*}{Lemma 2.6}
\newtheorem*{cor2.7*}{Corollary 2.7}
\newtheorem*{cor2.8*}{Corollary 2.8}
\newtheorem*{thm2.9*}{Theorem 2.9}
\newtheorem*{thm2.10*}{Theorem 2.10}
\newtheorem*{rem2.11*}{Remark 2.11}
\newtheorem*{prop2.12*}{Proposition 2.12}
\newtheorem*{def3.1*}{Definition 3.1}
\newtheorem*{lem3.2*}{Lemma 3.2}
\newtheorem*{cor3.3*}{Corollary 3.3}
\newtheorem*{lem3.4*}{Lemma 3.4}
\newtheorem*{cor3.5*}{Corollary 3.5}
\newtheorem*{thm3.6*}{Theorem 3.6}
\newtheorem*{thm3.7*}{Theorem 3.7}
\newtheorem*{lem3.8*}{Lemma 3.8}
\newtheorem*{lem3.9*}{Lemma 3.9}
\newtheorem*{prop3.10*}{Proposition 3.10}
\newtheorem*{thm3.11*}{Theorem 3.11}
\newtheorem*{lem4.1*}{Lemma 4.1}
\newtheorem*{rem*}{Remark}
\newtheorem*{cor4.2*}{Corollary 4.2}
\newtheorem*{thm4.3*}{Theorem 4.3}
\newtheorem*{lem4.4*}{Lemma 4.4}
\newtheorem*{lem4.5*}{Lemma 4.5}
\newtheorem*{thm4.6*}{Theorem 4.6}
\newtheorem*{lem4.7*}{Lemma 4.7}
\newtheorem*{cor4.8*}{Corollary 4.8}
\newtheorem*{thm4.9*}{Theorem 4.9}
\newtheorem*{lem4.10*}{Lemma 4.10}
\newtheorem*{thm4.11*}{Theorem 4.11}
\newtheorem*{def5.1*}{Definition 5.1}
\newtheorem*{lem5.2*}{Lemma 5.2}
\newtheorem*{lem5.3*}{Lemma 5.3}
\newtheorem*{cor5.4*}{Corollary 5.4}
\newtheorem*{def5.5*}{Definition 5.5}
\newtheorem*{rems*}{Remarks}
\newtheorem*{lem5.6*}{Lemma 5.6}
\newtheorem*{lem5.7*}{Lemma 5.7}
\newtheorem*{lem5.8*}{Lemma 5.8}
\newtheorem*{prop5.9*}{Proposition 5.9}
\newtheorem*{cor5.10*}{Corollary 5.10}
\newtheorem*{lem5.11*}{Lemma 5.11}
\newtheorem*{thm5.12*}{Theorem 5.12}
\newtheorem*{cor5.13*}{Corollary 5.13}
\newtheorem*{cor5.14*}{Corollary 5.14}
\newtheorem*{cor5.15*}{Corollary 5.15}
\newtheorem*{cor5.16*}{Corollary 5.16}
\newtheorem*{prop5.17*}{Proposition 5.17}
\newtheorem*{cor5.18*}{Corollary 5.18}
\newtheorem*{cor5.19*}{Corollary 5.19}
\newtheorem*{cor5.20*}{Corollary 5.20}
\newtheorem*{lem5.21*}{Lemma 5.21}
\newtheorem*{lem5.22*}{Lemma 5.22}
\newtheorem*{def5.23*}{Definition 5.23}
\newtheorem*{cor5.24*}{Corollary 5.24}
\newtheorem*{lem5.25*}{Lemma 5.25}
\newtheorem*{cor5.26*}{Corollary 5.26}

\numberwithin{equation}{section}

\newcommand{\BB}{\mathscr B}
\newcommand{\PP}{\mathscr P}
\newcommand{\SSS}{\mathscr S}
\newcommand{\TT}{\mathscr T}
\newcommand{\LL}{\mathcal L}
\newcommand{\RR}{\mathscr R}

\newcommand{\N}{\mathbb N}
\newcommand{\Q}{\mathbb Q}
\newcommand{\Z}{\mathbb Z}
\newcommand{\C}{\mathbb C}
\newcommand{\R}{\mathbb R}
\newcommand{\DD}{\mathscr D}

\title{Factors of type III and the distribution of prime numbers}

\author{Florin P. Boca}

\address[F. P. Boca]{Department of Mathematics, University of Wales Swansea, Swansea SA28PP}

\curraddr{School of Mathematics, Cardiff University, Senghennydd Road, Cardiff CF2 4YH}
\email{BocaFP@cardiff.ac.uk}

\address{Institute of Mathematics of the Romanian Academy, P. O. Box 1-764, Bucharest 70700, Romania (on leave)}

\author{Alexandru Zaharescu}

\address[A. Zaharescu]{Department of Mathematics, Massachusetts Institute of Technology, Cambridge, MA 02139, USA}

\curraddr{Institute for Advanced Study, School of Mathematics, Olden Lane, Princeton, NJ 08540, USA}
\email{zaharesc@math.ias.edu}

\address{Institute of Mathematics of the Romanian Academy, P. O. Box 1-764, Bucharest 70700, Romania (on leave)}

\date{November 1998}

\begin{document}

\maketitle

\section{Introduction}
In this paper we consider ITPFI factors (infinite tensor products of factors of type I)
\begin{equation*}
M_{\beta,\SSS}=\bigotimes\limits_{p\in\SSS}\  (\BB (\ell^2 (\N)),\phi (\beta,p)),
\end{equation*}
where $\beta \in (0,1]$, $\SSS$ is a subset of the set $\PP$ of all prime numbers and $\phi (\beta,p)$
is the normal state on $\BB(\ell^2(\N))$ with list of (multiplicity $1$) eigenvalues
$\lambda_{\beta,p,k}=(1-p^{-\beta})p^{-\beta k}$ for $k\in \N=\{ 0,1,2,\ldots\}$.
We denote $M_{\SSS}=M_{1,\SSS}$.

Our aim is to study the invariants S and T of these factors (see \cite[\S V]{Co4} for definitions).
Their isomorphism classes are related to the distribution of the elements of $\SSS$ in the presence of
the parameter $\beta \in (0,1]$. These kinds of factors can be viewed as `localizations' of $M_{\beta,\PP}$
which appear in the work of C. C. Moore \cite{Mo}, B. Blackadar \cite{Bl} for $\beta =1$, and J. B. Bost and
A. Connes \cite{BC,Co4} for $\beta >0$.

Let $\TT$ be the Toeplitz $C^*$-algebra generated by the unilateral shift $S$ acting on $\ell^2(\N)$ and $p>1$.
Consider the one-parameter group of automorphisms of $\TT$ defined by $\sigma_{t,p}(S)=p^{it}S$, with $t\in \R$.
Then, for any $\beta >0$, there exists a unique KMS$_\beta$-state $\phi (\beta,p)=\operatorname{Tr}(\cdot h_{\beta,p})$,
where $h_{\beta,p}$ denotes the trace class operator on $\ell^2(\N)$ with simple eigenvalues
$(1-p^{-\beta})p^{-\beta k}$ with $k\in \N$ (see \cite{BC}).
Since $\TT$ contains the compact operators on $\ell^2(\N)$, we have
$\pi_{\phi (\beta,p)} (\TT)^{\prime \prime} =\pi_{\phi(\beta,p)} (\BB(\ell^2(\N)))^{\prime\prime}$.
The ITPFI factors $M_{\beta,\PP}$ arise in \cite{BC} (see also \cite{Co4}) in connection with the
Hecke $C^*$-algebra $A=C^*_r (\Gamma_0 \backslash \Gamma \slash \Gamma_0)$ of the discrete group
\begin{equation*}
\Gamma =\bigg\{ \bigg( \begin{matrix} 1 & a \\ 0 & b \end{matrix} \bigg) :a\in \Q,b\in \Q_+^*\bigg\}
\end{equation*}
with respect to its almost normal subgroup
\begin{equation*}
\Gamma_0 =\bigg\{ \bigg( \begin{matrix} 1 & a \\ 0 & 1 \end{matrix} \bigg) : a\in \Z \bigg\}.
\end{equation*}

If $L(\gamma)$ and $R(\gamma)$ denote the cardinals of the image of the double coset $\Gamma_0 \gamma \Gamma_0$
in $\Gamma \slash \Gamma_0$ and in $\Gamma_0 \backslash \Gamma$, respectively, then
\begin{equation*}
\sigma_t (f)(\gamma) =\bigg( \frac{L(\gamma)}{R(\gamma)}\bigg)^{-it} f(\gamma) ,\quad
\mbox{\rm where $f\in A$, $\gamma\in\Gamma$, $t\in \R$,}
\end{equation*}
defines a $1$-parameter group of automorphisms of $A$. The $C^*$-algebra $A$ is generated by
\begin{equation*}
\mu_n =n^{-1/2} \chi \bigg( \Gamma_0 \bigg( \begin{matrix} 1 & 0 \\ 0 & n \end{matrix}\bigg) \Gamma_0 \bigg),
\quad \mbox{\rm for $n\in \N^* =\{ 1,2,3,\ldots \}$,}
\end{equation*}
and
\begin{equation*}
e(r)=  \chi \bigg( \Gamma_0 \bigg( \begin{matrix} 1 & r \\ 0 & 1 \end{matrix}\bigg) \Gamma_0 \bigg),
\quad \mbox{\rm for $r\in \Q /\Z$.}
\end{equation*}

It turns out that  the $C^*$-subalgebra $B$ generated by $\mu_n$, with $n\in \N^*$, is isomorphic to
$\bigotimes_{p\in \PP} \TT_p$, where $\TT_p=\TT$ for all $p\in \PP$, and $\sigma_t(\mu_n) =n^{it}\mu_n$.
Moreover, if $\beta >0$ and $\phi_\beta$ is a KMS$_\beta$-state of $(A, (\sigma_t)_{t\in\R})$, then
$\phi_\beta \vert_B$ coincides with $\bigotimes_{p\in\PP} \phi (\beta,p)$.

If $\Q_p$ denotes as usual the $p$-adic field and $\Z_p$ the subring of $p$-adic integers, then
\begin{equation*}
P(\Z_p) =\bigg\{ \bigg( \begin{matrix} 1 & a \\ 0 & b \end{matrix} \bigg) : a\in\Z_p, b\in\Z_p^\times \bigg\}
\end{equation*}
is a compact open subgroup of
\begin{equation*}
P(\Q_p)=\bigg\{ \bigg( \begin{matrix} 1 & a \\ 0 & b \end{matrix}\bigg) : a\in \Q_p, b\in \Q_p^\times \bigg\} .
\end{equation*}

By \cite[\S 3]{Bl} we have
\begin{equation*}
\operatorname{Sp}(\chi_{P(\Z_p)} \vert \LL (P(\Q_p)))=\{ (1-p^{-1})p^{-k}:k\in \N\}.
\end{equation*}

Furthermore, if $\SSS$ is a subset of $\PP$, the von Neumann algebra $\RR(P(\SSS))$ generated by the
(right) regular representation of the restricted direct product group
$P(\SSS)=\prod_{p\in\SSS} (P(\Q_p),P(\Z_p))$ is isomorphic with the ITPFI factor $M_\SSS$ (see also
Proposition 10 in \cite{BC}).

In \S 2 we study the case $\lambda \in (0,1]$ and prove the following.

\begin{thm1.1*}
\emph{For any $\beta,\lambda \in (0,1]$, there exists a subset $\SSS=\SSS_{\beta,\lambda}$ of $\PP$ such that
$M_{\beta,\SSS}$ is isomorphic with the Powers factor $R_\lambda$ if $\lambda \in (0,1)$ and with the
Araki-Woods type {\rm III}$_1$ factor $R_\infty$ if $\lambda=1$. Moreover, $\SSS$ can be chosen to have
density $0$ in $\PP$ and $M_{\beta,\PP}$ is isomorphic with $R_\infty$ for any $\beta \in (0,1]$.}
\end{thm1.1*}

A proof of the fact that the asymptotic ratio of $M_\PP$ coincides with  $[0,\infty)$ has been presented
in \cite[\S 4]{Bl}. Nevertheless, it is not clear to these authors why the different $q$ which appear in the
definition of the sets $I_n$  \cite[Theorem 4.4]{Bl} can be chosen to be distinct, as required in the definition
of the asymptotic ratio.

We call a subgroup $G$ of $\R$ $\beta$-representable if it is the T-group of an ITPFI factor of type
$M_{\beta,\SSS}$, where $\beta \in (0,1]$ and $\SSS$ is a subset of $\PP$. In \S 3, we prove that $1$-representable groups
are $\beta$-representable for any $\beta \in (0,1)$. Moreover, if $\beta_0>0.07$ and $0<\beta <\beta_0$, then any
$\beta_0$-representable subgroup is $\beta$-representable.

In \S 4 the class of $1$-representable subgroups of $\R$ is shown to be invariant under homotheties.
It is also proved that if $\beta \in (c_0,1)$, $G$ is $\beta$-representable and
$(1-\beta)/(1-c_0)<\lambda <\beta/c_0$, then $\lambda G$ is $\beta$-representable, where $c_0=0.535$.
The T-groups of the factors $M_{\beta,\SSS}$ for $\SSS \subset \PP$ and $\beta \in (0,1]$, are also shown to
coincide with the T-groups of certain ITPFI$_2$ factors (infinite tensor products of $M_2(\C)$) with largest eigenvalue
tending to $1$ and satisfying some explicit growth conditions (cf. Corollaries 3.2 and 3.8). For any
$\beta \in (0.5,1]$ and any $\SSS \subset \PP$, the factor $M_{\beta,\SSS}$ is shown to be ITPFI$_2$.

In the last section we show the existence of a large class of $1$-representable groups by proving the following result.

\begin{thm1.2*}
\emph{Let $\beta \in (0,1]$, let $G$ be a countable subgroup of $\R$, and let $\Sigma$
be a countable subset of $\R \setminus G$. Then, there exists a subset $\SSS=\SSS_{\beta,G,\Sigma}$ of $\PP$ such
that $\mbox{\rm T}(M_{\beta,\SSS})$ contains $G$ and does not intersect $\Sigma$.}
\end{thm1.2*}

This result should be compared with a corollary of the work of Connes and Woods  \cite{CW} on the characterization of
the flow of weights of ITPFI factors. It is a consequence of \cite{CW} (cf. \cite[Theorem 1.4]{GS}) that any subgroup
of $\R$ is the T-group of an ITPFI$_2$ factor.

From the previous discussion on the group von Neumann algebras $\RR(P(\SSS))$, we obtain the following corollary,
which answers the questions raised in \cite[p.271]{Bl}.

\begin{cor1.3*}
\emph{For any $\lambda \in [0,1]$, there exist a subset $\SSS=\SSS_\lambda$ of $\PP$ such that the group von Neumann
algebra $\RR(P(\SSS))$ is isomorphic with a hyperfinite type {\rm III}$_\lambda$ factor.}

\emph{Moreover, if $G$ is a countable subgroup of $\R$ and $\Sigma$ is a countable subset of $\R \setminus G=G^c$, then there exists
$\SSS=\SSS_{G,\Sigma}$, a subset of $\PP$, such that $G\subset \mbox{\rm T}(\RR (P(\SSS)))$ and
$\mbox{\rm T}(\RR(P(\SSS)))\cap \Sigma =\emptyset$.}
\end{cor1.3*}

By a result of J. Dixmier and L. Puk\'anszky \cite{Pu}, for any connected, separable, locally compact group $G$, the central
desintegration of $\RR(G)$ contains no type III factors. Connes' result on the uniqueness of the hyperfinite II$_\infty$ factor
\cite[Corollary 5]{Co2} shows that the only factors which can appear are either type I or isomorphic with $R_{0,1}$, the
hypefinite Araki-Woods type II$_\infty$ factor.

The first example of a locally compact group whose regular representation contains type III factors is due to
R. Godement: if $G=\R^2 \rtimes \operatorname{GL}_2(\Q)$, then $\RR(G)$ is a non-hyperfinite type III$_1$ factor.
As pointed out by Sutherland \cite{Su}, one can derive from this $G$, for any $\lambda \in(0,1)$, examples of
groups $G_\lambda$ with $\RR(G_\lambda)$ a non-hyperfinite type III$_\lambda$ factor and $G_{0,\lambda}$ such that
$\RR(G_{0,\lambda})$ is non-hyperfinite type III$_0$ and $\operatorname{T}(\RR(G_{0,\lambda}))=2\pi\Z/\log \lambda$.
Examples of groups $G$, $G_\lambda$ and $G_{0,\lambda}$, for $\lambda \in \{0,1\}$, with the same properties but producing
hyperfinite factors,were constructed by Connes (see the last part of \cite{Su}). In this respect, the previous corollary produces
the first examples of locally compact groups $G$ such that $\RR(G)$ is a hyperfinite type III$_0$ factor and
$\operatorname{T}(\RR(G))$ is dense in $\R$.

Another consequence of the separation theorem, Theorem 1.2, is that for any real number field $K$, there exists a $1$-representable
group $G\neq \{ 0\}$, invariant under homotheties with elements in the group of units of $K$.

{\bf Notation.} We use the following conventions:
\begin{itemize}
\item[]
$\vert F\rvert$ is the cardinality of a finite set $F$;
\item[]
$\pi(x)$ is the number of primes which do not exceed $x$;
\item[]
$\{ x\}$ is the fractional part of $x$;
\item[]
$\|x\|$ is the distance from $x$ to $\Z$;
\item[]
$f(x)=o(g(x))$ if $\lim_{x\rightarrow \infty}(f(x)/g(x))=0$;
\item[]
$f(x)=O(g(x))$ or  $f\ll g$ if $\lvert f(x)\rvert /g(x)$ is bounded at $+\infty$;
\item[]
$f(x)\sim g(x)$ if $\lim_{x\rightarrow \infty} (f(x)/g(x))=1$;
\item[]
The symbols $o$, $O$ and $\sim$ keep the same meaning for sequences;
\item[]
For two series with non-negative terms, we write
\begin{equation*}
\sum\limits_n a_n\sim\sum\limits_n b_n
\end{equation*}
if they are simultaneously divergent or convergent;
\item[]
If $X$ is a subset of $\R$, we denote $X^c=\R \setminus X =\{a\in\R: a\notin X\}$;
\item[]
If $B\subset \R^+$ and $\alpha \in \R$, $B^\alpha$ denotes the set $\{b^\alpha: b\in B\}$;
\item[]
For $\SSS$ a countable subset of $(0,\infty)$, we denote
\begin{equation*}
f_\SSS (t)=\sum\limits_{p\in\SSS} \frac{\sin^2 (2\pi t\log p)}{p} \in (0,\infty],\quad
\mbox{\rm for $t\in\R$.}
\end{equation*}
\end{itemize}

\newpage

\section{The type III$_{\lambda}$ case, $\lambda \in (0,1]$}

In this section we shall construct for any $\beta \in (0,1]$ subsets of primes which produce ITPFI
type III$_\lambda$ factors, $\lambda\in (0,1]$, of the type described in the introduction.

When $\beta=1$, we make use of the prime number theorem in the form (see for example, \cite[p.251]{Ve})
\begin{equation}\label{eq2.1}
\pi (x)=x/\log x+O(x/\log^2 x).
\end{equation}

When $\beta \in (0,1)$, we use results on the distribution of primes in intervals $[x,x+h]$,
where $h\in [x^c,x]$ and $c>0$ is a constant, of the following type: there exist $c_1,c_2,x_1>0$ such that
\begin{equation}\label{eq2.2}
\pi(x+h)-\pi(x) <c_1 h/\log x \quad \mbox{\rm for all $h\in [x^c,x]$, $x>x_1$},
\end{equation}
\begin{equation}\label{eq2.3}
\pi(x+h)-\pi(x) >c_2 h/\log x \quad \mbox{\rm for all $h\in [x^c,x]$, $x>x_1$}.
\end{equation}

Estimate \eqref{eq2.2} holds for any $c>0$, as a consequence of $\pi(x+h)-\pi(x) \ll h/\log x$, a result
proved by Hardy and Littlewood \cite{HL} using Brun's sieve. An important amount of work has been done to
lower the constant $c$ in \eqref{eq2.3}. The most recent results, due to R. C. Baker and G. Harman \cite{BH}, show
that $c$ can be lowered to $c_0=0.535$, whence there exists $x_1 >0$ such that
\begin{equation*}
\pi(x)-\pi(x-y) >y/20\log x \quad \mbox{\rm for all $y\geqslant x^{0.535}$, $x>x_1$}.
\end{equation*}

The existence of $c<\frac{2}{3}$ such that \eqref{eq2.3} is fulfilled (proved by Ingham in 1937) will suffice for the
purpose of this section.

We start with the following mere consequence of \eqref{eq2.1}.

\begin{lem2.1*}
\emph{Let $\{\varepsilon_n\}_n$ be a sequence of numbers in $(0,1)$ such that $\lim_n \varepsilon_n =0$ and
$\lim_n n\varepsilon_n=\infty$. Then}
\begin{equation*}
\lim\limits_n \frac{\pi(e^{(n+a+\varepsilon_n)/t})-\pi(e^{(n+a)/t})}{e^{n/t} \varepsilon_n/n}=e^{a/t},
\end{equation*}
\emph{for all $a\in \R$ and $t>0$.}
\end{lem2.1*}

\begin{proof}
Denote by $\ell$ the left-hand side limit. By \eqref{eq2.1} and the fact that
\begin{equation*}
\lim\limits_n \frac{nt^2 e^{(n+a+\varepsilon_n)/t}}{\varepsilon_n e^{n/t}(n+a+\varepsilon_n)^2}
=\lim\limits_n \frac{nt^2 e^{(n+a)/t}}{\varepsilon_n e^{n/t} (n+a)^2} =0,
\end{equation*}
it follows that
\begin{equation*}
\begin{split}
\ell & = \lim\limits_n \frac{nt}{\varepsilon_n e^{n/t}}
\bigg( \frac{e^{(n+a+\varepsilon_n)/t}}{n+a+\varepsilon_n} -\frac{e^{(n+a)/t}}{n+a} \bigg)
\\ & =te^{a/t} \lim\limits_n \frac{n((n+a)e^{\varepsilon_n /t}-(n+a+\varepsilon_n))}{\varepsilon_n (n+a+\varepsilon_n)(n+a)} \\
& =te^{a/t}\lim\limits_n \bigg( \frac{e^{\varepsilon_n/t}-1}{\varepsilon_n}-\frac{1}{n+a}\bigg) \\
& = te^{a/t} \lim\limits_n \frac{e^{\varepsilon_n/t}-1}{\varepsilon_n}=e^{a/t}. \qedhere
\end{split}
\end{equation*}
\end{proof}

\begin{cor2.2*}
\emph{Assume that $\{ \varepsilon_n\}_n$ is a sequence as in Lemma {\rm 2.1}, $a\in \R$, and $t_0 >0$, and take}
\begin{equation}\label{eq2.4}
B_n=\PP \cap ( e^{(n+a)/t_0},e^{(n+a+\varepsilon_n)/t_0}].
\end{equation}
\emph{Then}
\begin{equation*}
0< \inf\limits_n \frac{n}{\varepsilon_n} \sum\limits_{p\in B_n} \frac{1}{p} \leqslant
\sup\limits_n \frac{n}{\varepsilon_n} \sum\limits_{p\in B_n} \frac{1}{p} < \infty .
\end{equation*}
\emph{In particular, if we also assume that $\sum_n \varepsilon_n/n=\infty$ and take $B=\bigcup_n B_n$, then}
\begin{equation*}
\sum\limits_{p\in B} \frac{1}{p} \sim \sum\limits_n \frac{\varepsilon_n}{n}=\infty.
\end{equation*}
\emph{If in addition $\sum_n \varepsilon_n^3/n<\infty$ and $a\in \{ 0,\frac{1}{2}\}$, then}
\begin{equation*}
f_B (t_0) \sim \sum\limits_n \frac{\varepsilon_n^3}{n} <\infty .
\end{equation*}
\end{cor2.2*}

\begin{proof}
Lemma 2.1 yields
\begin{equation*}
\sum\limits_{p\in B_n} \frac{1}{p} \leqslant \frac{\lvert B_n\rvert}{e^{(n+a)/t_0}} =
\frac{\pi(e^{(n+a+\varepsilon_n)/t_0})-\pi( e^{(n+a)/t_0})}{e^{(n+a)/t_0}}
\leqslant \frac{c\varepsilon_n e^{n/t_0}}{ne^{(n+a)/t_0}}
\leqslant \frac{c_1\varepsilon_n}{n},
\end{equation*}
and similarly
\begin{equation*}
\sum\limits_{p\in B_n}\frac{1}{p} \geqslant \frac{c\varepsilon_n e^{n/t_0}}{ne^{(n+a+\varepsilon_n)/t_0}}
\geqslant \frac{c_2\varepsilon_n}{n},
\end{equation*}
for some constants $c,c_1,c_2 >0$ depending on $a$ and $t_0$.

When $a\in \{ 0,\frac{1}{2}\}$, we use $\lvert \sin (2\pi x)\rvert \leqslant 2\pi \operatorname{dist} (x,\frac{1}{2}\Z)$,
for $x\in \R$, to get
\begin{equation*}
f_B (t_0) \leqslant 4\pi^2 \sum\limits_n \varepsilon_n^2 \sum\limits_{p\in B_n} \frac{1}{p}
\sim\sum\limits_n \frac{\varepsilon_n^3}{n}.\qedhere
\end{equation*}
\end{proof}

To sort out the case $0<\beta <1$ one has to proceed more carefully. Firstly, fix $\beta \in(0,1)$, then fix
$\beta_0 \in [\beta,\frac{1}{3}(\beta+2)] \cap (c_0=0.535,1]$. Fix also $t_0>0$ and $a\in\R$ and consider
a sequence $\{ \varepsilon_n\}_n$ such that
\begin{equation}\label{eq2.5}
(\log n)^{-1} e^{-(1-\beta_0)n/\beta t_0} \leqslant \varepsilon_n \leqslant
(\log n)^{-1} e^{-(1-\beta)n/3\beta t_0}.
\end{equation}

Take $n$ large enough such that $\varepsilon_n <\beta t_0$ and define
\begin{equation*}
x_n=e^{(n+a)/t_0},\qquad h_n =e^{(n+a+\varepsilon_n)/\beta t_0} -e^{(n+a)/\beta t_0} \sim \varepsilon_n
e^{(n+a)/\beta t_0} /\beta t_0 ,
\end{equation*}
\begin{equation}\label{eq2.6}
B_n =\PP \cap (x_n,x_n+h_n] \qquad \mbox{\rm (mutually disjoint sets),}
\end{equation}
\begin{equation}\label{eq2.7}
B= \bigcup\limits_n B_n ,
\end{equation}
\begin{equation*}
a_n= \frac{\pi (x_n+h_n)-\pi(x_n)}{h_n/\log x_n}
=\frac{\lvert B_n\rvert}{h_n/\log x_n} ,\qquad
b_n =\frac{\lvert B_n\rvert}{\varepsilon_n e^{n/\beta t_0} /n} .
\end{equation*}

Since $\beta_0 >c_0$, we have $x_n \geqslant h_n \geqslant x_n^{c_0}$ for $n$ larger than some $n_1$. Therefore
by \eqref{eq2.2} and \eqref{eq2.3} the sequences $a_n$ and $1/a_n$ are bounded and so are $b_n$ and $1/b_n$.
This shows that, for some constant $c>0$,
\begin{equation}\label{eq2.8}
\sum\limits_{p\in B} \frac{1}{p^\beta} =\sum\limits_n \sum\limits_{p\in B_n} \frac{1}{p^\beta} \sim
\sum\limits_n \frac{\lvert B_n\rvert}{e^{n/t_0}} \geqslant c \sum\limits_n \frac{\varepsilon_n e^{n/\beta t_0}}{ne^{n/t_0}}
\geqslant c\sum\limits_n \frac{1}{n\log n} =\infty .
\end{equation}

When $a\in \{ 0,\frac{1}{2}\}$, we obtain, as at the end of the proof of Corollary 2.2,
\begin{equation}\label{eq2.9}
\begin{split}
f_{B^\beta}(t_0) = \sum\limits_{p\in B} \frac{\sin^2 (2\pi \beta t_0 \log p)}{p^\beta} &
\leqslant 4\pi^2 \sum\limits_n \sum\limits_{p\in B_n} \frac{\varepsilon_n^2}{p^\beta}
\sim 4\pi^2 \sum\limits_n \frac{\varepsilon_n^2 \lvert B_n\rvert}{e^{n/t_0}} \\
& \sim \sum\limits_n \frac{\varepsilon_n^3 e^{n/\beta t_0}}{ne^{n/t_0}}
\leqslant \sum\limits_n \frac{1}{n\log^3 n} <\infty ,
\end{split}
\end{equation}
whence we get the following proposition.

\begin{prop2.3*}
\emph{Let $\beta\in (0,1)$, $t_0>0$, $a\in\R$ and let $\{ \varepsilon_n\}_n$ be a sequence as in \eqref{eq2.5}.
Define $B_n$ and $B$ as in \eqref{eq2.6} and \eqref{eq2.7}. Then}
\begin{equation*}
\sum\limits_{p\in B} \frac{1}{p^\beta} =\infty .
\end{equation*}
\emph{Moreover, if we also assume $a\in \{ 0,\frac{1}{2}\}$, then}
\begin{equation*}
f_{B^\beta} (t_0) =\sum\limits_{p\in B} \frac{\sin^2 (2\pi \beta t_0 \log p)}{p^\beta} < \infty .
\end{equation*}
\end{prop2.3*}

Notice that we can actually use any $c_0<\frac{1}{3}$ instead of $c_0=0.535$.

\begin{rem2.4*}
Let $\beta \in (0,1]$ and $t_0>0$, $a$, $\{ \varepsilon_n\}_n$, $B_n$ and $B$ be as in Proposition 2.3.
Then, if $\SSS$ is a subset of $\N^*$ of positive density, we have
\begin{equation*}
\sum\limits_{n\in\SSS} \sum\limits_{p\in B_n}\frac{1}{p^\beta}=\infty .
\end{equation*}
\end{rem2.4*}

\begin{prop2.5*}
\emph{Let $\beta \in (0,1]$ and $t_0>0$. Then, there exists a subset $B=B(\beta,t_0)$ of $\PP$ such that }
\begin{equation*}
f_{B^\beta} (t) <\infty \  \ \Longleftrightarrow \ \  t\in\Z t_0 .
\end{equation*}
\end{prop2.5*}
.
\begin{proof}
Take $\varepsilon_n=(\log n)^{-1}e^{-(1-\beta)n/3\beta t_0}$ and $n$ such that $\varepsilon_n <\beta t_0$.
The sets of primes defined by $B_n=\PP \cap [e^{(n+1/2)/\beta t_0},e^{(n+1/2+\varepsilon_n)/\beta t_0})$ are
mutually disjoint. We set $B=\bigcup_n B_n$.

Since $\frac{1}{2} <\{ \beta t_0\log p\} <\varepsilon_n+\frac{1}{2}$, for $p\in B_n$, we get
$\sin^2 (2\pi \beta t_0 \log p) \leqslant 4\pi^2 \varepsilon_n^2$. By Proposition 2.3,
\begin{equation*}
\sum\limits_{p\in B} \frac{1}{p^\beta} =\infty \quad \mbox{\rm and} \quad
f_{B^\beta} (t_0)=\sum\limits_{p\in B} \frac{\sin^2 (2\pi\beta t_0\log p}{p^\beta } <\infty .
\end{equation*}

Since $\sin^2(x+y) \leqslant 2(\sin^2 x +\sin^2 y)$, the set $G_{B^\beta} =\{ t\in\R: f_{B^\beta} (t) <\infty\}$
is s subgroup of $\R$ and $\Z t_0 \subset G_{B^\beta}$.

Finally, we show that if $t\notin \Z t_0$, then $f_{B^\beta} (t)=\infty$. Let $t=mt_0+t^\prime \notin \Z t_0$ with
$m\in\Z$ and $0<t^\prime<t_0$. Since $mt_0\in G_{B^\beta}$, then $t\in G_{B^\beta}$ if and only if $t^\prime \in G_{B^\beta}$.
Denote $\rho=t^\prime/t_0 \in(0,1)$. Since $n+\frac{1}{2} <\beta t_0 \log p < \varepsilon_n +n+\frac{1}{2}$, for $p\in B_n$,
we get, for all $p\in B_n$,
\begin{equation*}
\pi\rho +2\pi n\rho <2\pi \beta t^\prime \log p <\pi\rho+2\pi\rho \varepsilon_n +2\pi n\rho.
\end{equation*}

Using Weyl's uniform distribution theorem \cite{We} for irrational $\rho$, we clearly see that,
in both cases (irrational or rational $\rho$), the set
\begin{equation*}
S_\rho =\{ n\in\N^*: \tfrac{2}{3} \pi \geqslant \{ 2\pi \beta t^\prime \log p\}\geqslant \tfrac{1}{3} \pi \
\mbox{\rm for all $p\in B_n$}\}
\end{equation*}
has positive density in $\N^*$. Using Remark 2.4 also we get
\begin{equation*}
f_{B^\beta}(t^\prime) \geqslant \sum\limits_{n\in S_\rho} \sum\limits_{p\in B_n}
\frac{\sin^2 (2\pi\beta t^\prime \log p)}{p^\beta} \geqslant
\sum\limits_{n\in S_\rho} \sum\limits_{p\in B_n} \frac{3}{4p^\beta}=\infty .
\end{equation*}

Remark that we have used $n+\frac{1}{2}$ rather than $n$ in order to cover the case $\rho=\frac{1}{2}$ as well.
\end{proof}

The next lemma shows a connection between $\operatorname{T}(M_{\beta,\SSS})$ and the distribution
of $\log p$, for $p\in\SSS$.

\begin{lem2.6*}
\emph{Let $\SSS$ be a subset of $\PP$ and $\beta \in (0,1]$. Then}
\begin{equation*}
T(M_{\beta,\SSS})=\bigg\{ t\in \R: \sum\limits_{p\in\SSS} \frac{1}{p^\beta} \sin^2 \bigg(\frac{\beta t\log p}{2}\bigg) <\infty \bigg\} .
\end{equation*}
\end{lem2.6*}

\begin{proof}
Using \cite[\S V.4]{Co4} one can easily check that
\begin{equation*}
M_{\beta,\SSS} \ \mbox{\rm is of type III}\ \  \Longleftrightarrow \ \
\sum\limits_{p\in\SSS} \frac{1}{p^\beta}=\infty ,
\end{equation*}
in which case $\operatorname{T}(M_{\beta,\SSS}) \neq \R$ (cf. \cite{Ta}).
By \cite[Corollary 1.3.9]{Co1} we know that
\begin{equation*}
t\in\operatorname{T}(M_{\beta,\SSS}) \\  \Longleftrightarrow \ \
\sum\limits_{p\in\SSS}\bigg( 1-\bigg| \sum\limits_{k\geqslant 0} (p^{-\beta k} (1-p^{-\beta}))^{1+it}\bigg|\bigg) <\infty .
\end{equation*}
Therefore
\begin{equation}\label{eq2.10}
t\in\operatorname{T}(M_{\beta,\SSS}) \ \  \Longleftrightarrow \ \
\sum\limits_{p\in\SSS} \bigg( 1-\frac{1-p^{-\beta}}{\lvert 1-p^{-\beta (1+it)}\rvert}\bigg) < \infty .
\end{equation}

As $\lim_p \lvert 1-p^{-\beta (1+it)} \rvert =1$, the right-hand side of \eqref{eq2.10} converges if and only if
the series $A_\SSS =\sum_{p\in\SSS} (\lvert 1-p^{-\beta (1+it)}\rvert -(1-p^{-\beta}))$ does.

If we set $x=p^{-\beta}$, $\theta=\theta(x)=-\beta t\log p$, then
\begin{equation*}
\lvert 1-p^{-\beta(1+it)}\rvert -(1-p^{-\beta})=\frac{4x\sin^2(\theta/2)}{1-x+\sqrt{1-2x\cos\theta+x^2}} .
\end{equation*}
As $\lim_{x\rightarrow 0}(1-x+\sqrt{1-2x\cos\theta+x^2})=2$, $A_\SSS$ converges if and only if
$\sum_{p\in\SSS}p^{-\beta}\sin^2 (\frac{1}{2}\beta t\log p)$ does.
\end{proof}

\begin{cor2.7*}
\emph{For all $\beta \in (0,1]$ and $\lambda\in(0,1)$, there exists a subset $B=B(\beta,\lambda)$ of $\PP$ such
that $\operatorname{T}(M_{\beta,B})=2\pi\Z /\log \lambda$.}
\end{cor2.7*}

\begin{proof}
Take $t_0$ such that $2t_0\log\lambda=-1$ and $B\subset \PP$ constructed in Proposition 2.5.
By Lemma 2.6, $\operatorname{T}(M_{\beta,B})=4\pi t_0 \Z =2\pi \Z/\log\lambda$.
\end{proof}

If we consider $M_{\beta,\PP}=\bigotimes_{p\in\PP} (\BB(\ell^2(\N)),\phi(\beta,p))$,the proof of Proposition 2.5 also shows that
$\sum_{p\in\PP} p^{-\beta} \sin^2(\beta t\log p)=\infty$ for all $t\neq 0$ (actually this can be established using only
$\pi (x) \sim x/\log x$ instead of \eqref{eq2.1}),which gives the following corollary.

\begin{cor2.8*}
\emph{The equality $\operatorname{T}(M_{\beta,\PP})=\{ 0\}$ holds for all $\beta \in (0,1]$.}
\end{cor2.8*}

This kind of argument will enable us to prove that for any $\beta \in (0,1]$, the ITPFI factor
$M_{\beta,B}$, with $B=B(\beta,\lambda)$ as in Corollary 2.7, is isomorphic to the Powers factor $R_\lambda$
for each $\lambda\in (0,1)$ and also for any $\beta\in (0,1]$, $M_{\beta,\PP}$ is isomorphic to the Araki-Woods
type III$_1$ factor $R_\infty$, known to be the unique hyperfinite factor of type III$_1$ by the deep
work of A. Connes and U. Haagerup \cite{Co3,Ha}. By Connes' duality between the S and T invariants
\cite[Theorem 3.4.1]{Co1}, we only have to construct elements in the asymptotic ratio set
which are different from $0$ and $1$.

\begin{thm2.9*}
\emph{Let $\lambda\in (0,1)$, $\beta \in (0,1]$ and $B=B(\beta,\lambda)$ be the subset of $\PP$ from Corollary 2.7. Then}
\begin{equation*}
S(M_{\beta,B})=\{ 0\} \cup \lambda^\Z .
\end{equation*}
\end{thm2.9*}

\begin{proof}
In order to prove that $\operatorname{S}(M_{\beta,B})=\{ 0\} \cup \lambda^\Z$, it suffices to show that
$\lambda^{-2} \in\operatorname{S}(M_{\beta,B})=r_\infty (M_{\beta,B})$, and then use \cite[Theorem 3.4.1]{Co1},
for we already know that $\operatorname{T}(M_{\beta,B})=2\pi\Z /\log \lambda$.

Let $t_0>0$ be such that $2t_0\log\lambda=-1$ and $\varepsilon_n=\varepsilon_n (\beta,t_0)$ be the sequence used
in the proof of Proposition 2.5. Consider also the (mutually disjoint) subsets $B_n$ from the proof of Proposition 2.5
and denote their union by $B$. From the definition of $B_n$ we see that
\begin{equation}\label{eq2.11}
e^{(1-\varepsilon_n)/t_0} < p^\beta/q^\beta < q^{(1+\varepsilon_n)/t_0},\quad
\mbox{\rm for $p\in B_{2n+1}$, $q\in B_{2n}$.}
\end{equation}

As in \eqref{eq2.8}, we find constants $c_1,c_2>0$ such that for all $n\in\N$, $n\geqslant 2$,
\begin{equation}\label{eq2.12}
c_1(n\log n)^{-1} e^{(\beta+2)n/3\beta t_0} <\lvert B_n \rvert < c_2(n\log n)^{-1}e^{(\beta+2)n/3\beta  t_0},
\end{equation}
whence $\lim_n\lvert B_n\rvert=\infty$ and $\inf_n \lvert B_{n+1}\rvert /\lvert B_n\rvert >0$.
We may subsequently select $N$ and nonempty subsets $B_n^\prime \subset B_n$, $n\geqslant N$,
such that $\inf_n \lvert B_n^\prime\rvert /\lvert B_n\rvert >0$ and $\lvert B_{2n}^\prime \rvert =
\lvert B_{2n+1}^\prime \rvert$ for all $n\geqslant N$. We consider
\begin{equation*}
B^\prime =\bigcup\limits_n B^\prime_{2n+1}=\{ p_m\}_m,\quad
B^{\prime\prime} =\bigcup\limits_n B^{\prime}_{2n}=\{ q_m\}_m \subset B.
\end{equation*}
The elements $p_j$, $q_k$ are all distinct and, from \eqref{eq2.11}, $\lim_m p_m^\beta/q_m^\beta =e^{1/t_0}=\lambda^{-2}$.
Moreover, \eqref{eq2.12} and $\lim_n\varepsilon_n =0$ provide a constant $c_3>0$ such that
\begin{equation*}
\begin{split}
\sum\limits_m \frac{1}{q_m^\beta} & = \sum\limits_{p\in B^{\prime \prime}} \frac{1}{p^\beta}
= \sum\limits_n \sum\limits_{p\in B_{2n}^\prime} \frac{1}{p^\beta} >
\sum\limits_n \frac{\lvert B_{2n}^\prime\rvert}{e^{(\varepsilon_{2n}+2n+1/2)/t_0}}
 >c_3 \sum\limits_n \frac{e^{2n(\beta+2)/3\beta t_0}}{n\log n} \\
 & =c_3\sum\limits_n \frac{e^{(4(1-\beta)n/3\beta t_0)}}{n\log n}
 \geqslant c_3 \sum\limits_n \frac{1}{n\log n}=\infty ,
\end{split}
\end{equation*}
which also entrains $\sum_m p_m^{-\beta}=\infty$. Summarizing, we have produced two sequences
of disjoint primes $\{ p_m\}_m$ and $\{q_m\}_m$ in $B$ such that
\begin{equation*}
\lim\limits_m \frac{p_m^\beta}{q_m^\beta}=\lambda^{-2} \quad \mbox{\rm and} \quad
\sum\limits_m\frac{1}{p_m^\beta}=\infty =\sum\limits_m \frac{1}{q_m^\beta} .
\end{equation*}

Taking $K_m^1$and $K_m^2$ to be the singleton sets containing $q_n^{-\beta}(1-p_n^{-\beta})(1-q_n^{-\beta})$
and $p_n^{-\beta}(1-p_n^{-\beta})(1-q_n^{-\beta})$, respectively, and $I_n=\{ p_n,q_n\}$, we conclude by the
definition of the asymptotic ratio set \cite{AW} that $\lambda^{-2} \in r_\infty (M_{\beta,B})$.
\end{proof}

\begin{thm2.10*}
\emph{For any $\beta \in (0,1]$ we have}
\begin{equation*}
\operatorname{S}(M_{\beta,\PP}) = [0,\infty).
\end{equation*}
\emph{In particular, all factors $M_{\beta,\PP}$ are isomorphic with $R_\infty$.}
\end{thm2.10*}

\begin{proof}
We have to prove that $(0,1] \subset r_\infty (M_{\beta,\PP})$.
Let $\lambda\in (0,1)$, $\varepsilon_n=(\log n)^{-1}$, $a=\log\lambda$ and
$B_{n,a}=\PP \cap (e^{(n+a)/\beta},e^{(n+a+\varepsilon_n)/\beta}]$, which are mutually disjoint and
\begin{equation}\label{eq2.13}
e^{(a-\varepsilon_n)/\beta} < p/q < e^{(a+\varepsilon_n)/\beta}, \qquad \mbox{\rm for $p\in B_{n,a}$, $q\in B_{n,0}$}.
\end{equation}

By Lemma 2.1, $\lim_n \lvert B_{n,a}\rvert /\lvert B_{n,0}\rvert =e^{a/\beta} <1$ and we can subsequently assume
$\lvert B_{n,a}\rvert \leqslant \lvert B_{n,0}\rvert$ for all $n$. We choose for each $n$ a subset
$B^\prime_{n,0} \subset B_{n,0}$ such that $\lvert B^\prime_{n,0}\rvert =\lvert B_{n,0}\rvert$, and then consider
the disjoint subsets $B^\prime =\bigcup_n B^\prime_{n,0}=\{ p_m\}_m$ and
$B^{\prime \prime} =\bigcup_n B_{n,a}=\{ q_m\}_m$ of $\PP$. By \eqref{eq2.13}, $\lim_m p_m^\beta /q_m^\beta =a^a =\lambda$.
By construction, $p_j$ and $q_k$ are all distinct. Moreover, since
\begin{equation*}
\sum\limits_m \frac{1}{q_m^\beta} \geqslant \sum\limits_m \frac{1}{q_m} =\sum\limits_n \sum\limits_{p\in B_{n,a}}
\frac{1}{p} \sim \sum\limits_n \frac{\varepsilon_n}{n} =\infty ,
\end{equation*}
it follows that $\sum_m p_m^{-\beta} =\infty$. If $K_n^1$, $K_n^2$ and $I_n$ are defined as in the proof of Theorem 2.9,
we conclude that $\lambda \in r_\infty (M_{\beta,\PP})$ by the definition of the asymptotic ratio.
\end{proof}

\begin{rem2.11*}
(i) Let $\lambda,\mu\in (0,1)$, with $\lambda \neq \mu$. We can extract two disjoint subsets
$B^\prime (\beta,\lambda) \subset B(\beta,\lambda)$ and $B^\prime (\beta,\mu) \subset B(\beta,\mu)$ such that
$M_{B^\prime (\beta,\lambda)}$ is isomorphic to $R_\lambda$ and $M_{B^\prime (\beta,\mu)}$ to $R_\mu$.
Since $R_\lambda \otimes R_\mu$ is isomorphic to $R_\infty$ if $\log\lambda /\log\mu$ is irrational,
and $R_\infty$ is isomorphic to $R_\infty \otimes M$ for any hyperfinite type III factor $M$, we obtain
yet another proof of the fact that $M_{\beta,\PP}=\bigotimes_{p\in \PP} (\BB(\ell^2(\N)),\phi (\beta,p))$ is
of type $\operatorname{III}_1$.

(ii) The set $B(\beta,\lambda)$ has density $0$ in $\PP$ for any $\lambda \in (0,1)$ and $\beta \in (0,1]$, for an
immediate computation shows that, for any $a\in \R$ and any $t_0 >0$,
\begin{equation*}
\lim\limits_n \frac{\sum_{k=1}^n (\pi(e^{(k+a+\varepsilon_k)/t_0}) -\pi(e^{(k+a)/t_0}))}{\pi(e^{(n+a+\varepsilon_n)/t_0})} =0.
\end{equation*}
\end{rem2.11*}

By a theorem of J. Dixmier and M. Takesaki \cite{Ta}, we know that if $\sum_{p\in\SSS} p^{-\beta} =\infty$ for some subset
$\SSS$ of $\PP$, then $\operatorname{T} (M_{\beta,\SSS}) \neq \R$ (for in this case $M$ is type III, and thus not hyperfinite).
We conclude this section by giving an elementary proof of this fact.

\begin{prop2.12*}
\emph{Let $\{ p_k\}_k$ and $\{ a_k\}_k$ be two sequences such that $p_k,a_k\geqslant 0$,
$\sum_k p_k^{-1} =\infty$ and $\lim_k a_k=\infty$. If $h:[0,1]\rightarrow \R_+^*$ is a function such that
$h(t) \geqslant \min (h(c),h(1-c)) >0$ for any $t\in [c,1-c]$ with $0\leqslant c<\frac{1}{2}$, then the set }
\begin{equation*}
S=\bigg\{ t\in \R: \sum\limits_k \frac{h(\{ ta_k\})}{p_k} < \infty \bigg\}
\end{equation*}
\emph{has Lebesgue measure $0$ in $\R$.}
\end{prop2.12*}

\begin{proof}
Denote
\begin{equation*}
\begin{split}
f_n(t) & = \sum\limits_{k=1}^n \frac{h(\{ ta_k\})}{p_k},\quad \mbox{\rm for $t\in\R$,} \\
f(t) & = \lim\limits_n f_n(t) \quad \mbox{\rm (which is finite on $S$),} \\
V_n & =\sum\limits_{k=1}^n \frac{1}{p_k} , \\
f_{n,c}(t) & =\sum\limits_{k=1}^n \frac{\chi_{[c,1-c]} (\{ ta_k\})}{p_k} ,
\quad \mbox{\rm for $0<c<\frac{1}{2}$, $t\in [0,1]$,} \\
m_c & = \min (h(c),h(1-c)) >0.
\end{split}
\end{equation*}

Clearly $V_n \geqslant f_{n,c}(t)$ and $V_n \geqslant f_n(t) \geqslant m_c f_{n,c}(t)$ for all $n,c,t$.
Let $\mu$ be the Lebesgue measure on $[0,1]$. For any $c\in (0,\frac{1}{2})$,
\begin{equation*}
\lim\limits_{x\rightarrow \infty} \mu (\{ t\in [0,1]: \{ tx\} \in [0,c]\cup [1-c,1)\}) =2c,
\end{equation*}
and since $\lim_k a_k =\infty$, there exists $n_c$ such that, for all $k\geqslant n_c$,
\begin{equation*}
\mu (\{ t\in [0,1]: \{ ta_k\} \in [0,c] \cup [1-c,1) \}) \leqslant 3c.
\end{equation*}
Therefore, for all $n\geqslant n_c$,
\begin{equation}\label{eq2.14}
\int_0^1 f_{n,c}(t)\, dt =V_n -\sum\limits_{k=1}^n \int_0^1 \frac{\chi_{[0,c]\cup [1-c,1)} (\{ ta_k\})}{p_k}\, dt
\geqslant V_n -V_{n_c} -3c(V_n-V_{n_c}) .
\end{equation}

If we set $A_{n,,c} =\{ t\in [0,1]: f_{n,c}(t) < \frac{1}{2} V_n\}$ and estimate the integral of $f_{n,c}$ on $A_{n,c}$
and on its complement, then
\begin{equation}\label{eq2.15}
\int_0^1 f_{n,c}(t)\, dt \leqslant \tfrac{1}{2} V_n \mu (A_{n,c}) +V_n (1-\mu(A_{n,c})).
\end{equation}

We compare \eqref{eq2.14} and \eqref{eq2.15} to get, for all $n\geqslant n_c$,
\begin{equation*}
1-\tfrac{1}{2} \mu(A_{n,c}) \geqslant (1-3c)(1-V_{n_c}/V_n);
\end{equation*}
thus there exists $N_c \geqslant n_c$ such that
\begin{equation*}
1-\tfrac{1}{2} \mu(A_{n,c}) \geqslant 1-4c,\quad \mbox{\rm for $n\geqslant N_c$,}
\end{equation*}
and consequently
\begin{equation*}
\mu(A_{n,c}) \leqslant 8c, \quad \mbox{\rm for $n\geqslant N_c$.}
\end{equation*}

To conclude, let $\varepsilon >0$ and pick $c(k)\in (0,\frac{1}{2})$ with $\sum_k c(k) < \frac{1}{8} \varepsilon$.
For each $k$ also choose $n(k) \geqslant N_{c(k)}$ such that $V_{n(k)} \geqslant 2k/m_{c(k)}$. Then
$\sum_k \mu (A_{n(k),c(k)})\leqslant \sum_k 8c(k) \leqslant \varepsilon$. Also, if $t_0 \in [0,1] \setminus
\bigcup_k A_{n(k),c(k)}$, then we get $f_{n(k),c(k)} (t_0) > \frac{1}{2} V_{n(k)}$ and
\begin{equation*}
f(t_0) \geqslant m_{c(k)} f_{n(k),c(k)} (t_0) \geqslant \tfrac{1}{2} m_{c(k)} V_{n(k)} > k,
\end{equation*}
for all $k$, showing that $S\cap [0,1] \subset \bigcup_k A_{n(k),c(k)}$, so $\mu (S \cap [0,1]) \leqslant \varepsilon$.
One shows in a similar way that $\mu (S \cap [n,n+1]) =0$ for all $n\in\Z$.
\end{proof}

The next argument \cite[Remark 1.2(b)]{GS} shows that $\operatorname{T}(M_{\beta,\SSS})$ is a $K_\sigma$ set for any
$\beta \in (0,1]$ and any set $\SSS \subset \PP$. For, consider the compact subsets
\begin{equation*}
K_{m,N} =\bigg\{ t\in \R: \lvert \theta \rvert +\sum\limits_{p\in\SSS, p\leqslant N} \frac{1}{p^\beta}
\sin^2 \bigg( \frac{\beta t\log p}{2}\bigg) \leqslant m \bigg\}
\end{equation*}
of $\R$. Then $K_m =\bigcap_N K_{m,N}$ is compact too and $\operatorname{T}(M_{\beta,\SSS})$ is $K_\sigma$.

\section{$\beta$-representable subgroups of $\R$}
We shall modify slightly the definition from the introduction and set, for any countable subset
$A\subset (0,\infty)$ and any $t\in\R$,
\begin{equation*}
f_A (2t) =\sum\limits_{n\in A} \frac{\sin^2 (t\log n)}{n}\in [0,\infty] .
\end{equation*}

The set $G_A=\{ t\in \R: f_A(t) < \infty\}$ is a subgroup of $\R$ since
\begin{equation}\label{eq3.1}
\sin^2 (x+y) \leqslant 2(\sin^2 x+\sin^2 y),\quad \mbox{\rm for $x,y\in \R$.}
\end{equation}

If $A\subset \PP^\beta=\{ p^\beta:p\in\PP\}$, then $G_A=\operatorname{T}(M_{\beta,A})$ by Lemma 2.6.

If $a=\{ a_k\}_k$ is a sequence in $\R^*_+$, then
\begin{equation*}
G(a)=\bigg\{ t\in\R: \sum\limits_k \frac{1}{a_k} \sin^2 \bigg( \frac{t\log a_k}{2}\bigg) <\infty\bigg\}
\end{equation*}
is also a subgroup of $\R$ by \eqref{eq3.1}.

\begin{def3.1*}
Let $\beta \in (0,1]$. S subgroup $G$ of $\R$ is called
\begin{itemize}
\item[(i)] \emph{$\beta$-representable} if there exists $A\subset \PP^\beta$ such that $G=G_A$.
\item[(ii)]
\emph{$\beta$-admissible} if it is the intersection of a family of $\beta$-representable groups.
\end{itemize}
Groups which are $1$-representable, or $1$-admissible, will simply be called \emph{representable} or
\emph{admissible}, respectively.
\end{def3.1*}

First, we prove that representable groups are also $\beta$-representable for any $0<\beta<1$.
The proof uses the fact that the function $\phi_t(x)=x\sin^2(t\log x)$, for $x\in (0,\infty)$,
is $C^1$ with $\lvert \phi_t^\prime (x)\rvert \leqslant 2t+1$ for all $x>0$. Changing $x$
into $x^{-1}$ we get, for all $x,y\in (0,\infty)$,
\begin{equation}\label{eq3.2}
\bigg| \frac{\sin^2(t\log x)}{x}-\frac{\sin^2(t\log y)}{y}\bigg| \leqslant
(2t+1)\bigg| \frac{1}{x}-\frac{1}{y}\bigg|,
\end{equation}
so if $a=\{ a_k\}_k$ and $b=\{ b_k\}_k$ are sequences from $\R_+^*$ and
$\sum_k \lvert a_k^{-1}-b_k^{-1}\rvert <\infty$, then
\begin{equation*}
\sum\limits_k \frac{\sin^2 (t\log a_k)}{a_k} \sim \sum\limits_k \frac{\sin^2 (t\log b_k)}{b_k}
\end{equation*}
for all $t\in\R$, and therefore $G(a)=G(b)$.

For each $c>0$ and $s>1$ consider the sequence $y=y(s,c)=\{ y_n\}_n$ defined by
\begin{equation}\label{eq3.3}
y_n=y_n(s,c)=e^{cn/\log^s n} .
\end{equation}

The following features of this sequence will be used in the proof of Theorem 3.6 and in \S 4.

\begin{lem3.2*}
For all $c>0$ and $s>1$, the following hold:
\begin{itemize}
\item[(i)]
$\displaystyle \quad \frac{y_{n+1}-y_n}{y_n} \sim \frac{c}{\log^s n}$;
\item[(ii)]
$\displaystyle \quad \lim\limits_n \frac{y_n/\log^2 y_n}{y_{n+1}/\log y_{n+1} -y_n/\log y_n}
=\lim\limits_n \frac{y_{n+1}/\log^2 y_{n+1}}{y_{n+1}/\log y_{n+1} -y_n/\log y_n} =0$.
\end{itemize}
\end{lem3.2*}

\begin{proof}
(i) Let $f(x)=cx/\log^s x$ and $g(x)=f(x+1)-f(x)$. Then, there exists $\xi_x \in (x,x+1)$ such that
\begin{equation*}
g(x)=f^\prime (\xi_x)=c\frac{\log\xi_x -s}{\log^{s+1} \xi_x} ,
\end{equation*}
and consequently $g(x)\sim c/\log^s x$. In particular, $\lim_{x\rightarrow\infty} g(x)=0$ and
\begin{equation*}
y_{n+1}/y_n -1 =e^{g(n)} -1 \sim g(n) \sim c/\log^s n .
\end{equation*}

(ii) Denote by $\ell$ the first limit. Since $\lim_n y_{n+1}/y_n=1$, we obtain
\begin{equation*}
\ell=\lim\limits_n \frac{y_n}{y_{n+1}\log y_n -y_n\log y_{n+1}} =
\lim\limits_{x\rightarrow\infty} \frac{1}{f(x)e^{g(x)}-f(x+1)} .
\end{equation*}
Therefore
\begin{equation*}
\begin{split}
\frac{1}{\ell} & =\lim\limits_{x\rightarrow\infty} (f(x)e^{g(x)} -f(x+1))
=\lim\limits_{x\rightarrow\infty} (f(x)(e^{g(x)}-1)-g(x)) \\
& =\lim\limits_{x\rightarrow\infty} f(x)(e^{g(x)}-1)=\lim\limits_{x\rightarrow\infty} f(x)g(x)
=\lim\limits_{x\rightarrow \infty} \frac{c^2 x}{\log^{2s} x}=\infty.
\end{split}
\end{equation*}
The second equality follows in a similar way.
\end{proof}

\begin{cor3.3*}
\emph{If $\alpha_n=\pi(y_{n+1})-\pi(y_n)$, then}
\begin{equation*}
\alpha_n \sim \frac{y_{n+1}}{\log y_{n+1}} -\frac{y_n}{\log y_n} .
\end{equation*}
\end{cor3.3*}

\begin{lem3.4*}
\begin{equation*}
\frac{y_{n+1}}{\log y_{n+1}} -\frac{y_n}{\log y_n} \sim
\frac{cy_n}{(\log y_n)(\log \log y_n)^s} .
\end{equation*}
\end{lem3.4*}

\begin{proof}
With the notation from the proof of Lemma 3.2 we get
\begin{equation*}
\begin{split}
L & = \lim\limits_{x\rightarrow\infty} \frac{e^{f(x+1)}/f(x+1)-e^{f(x)}/f(x)}{(ce^{f(x)})/(f(x)\log^s f(x))}
=\lim\limits_{x\rightarrow \infty} \frac{f(x)e^{g(x)}-f(x+1)}{cf(x)}\log^s x \\
& =\frac{1}{c} \lim\limits_{x\rightarrow \infty} g(x) \bigg( \frac{e^{g(x)}-1}{g(x)} -\frac{1}{f(x)}\bigg) \log^s x
=\frac{1}{c} \lim\limits_{x\rightarrow\infty} g(x)\log^s x =1 .\qedhere
\end{split}
\end{equation*}
\end{proof}

\begin{cor3.5*}
\emph{If $y_n=y_n(s,c)$ and $\alpha_n =\pi(y_{n+1})-\pi(y_n)$, for $c>0$ and $s>1$, then}
\begin{equation*}
\alpha_n=\alpha_n (s,c) \sim \frac{cy_n}{(\log y_n)(\log \log y_n)^s} .
\end{equation*}
\end{cor3.5*}

\begin{thm3.6*}
\emph{Representable subgroups of $\R$ are $\beta$-representable for any $\beta \in (0,1)$.}
\end{thm3.6*}

\begin{proof}
Let $G=G_B$ for some subset $B$ of $\PP$ and $\beta \in (0,1)$. We construct a subset $C$ of $\PP$ and a
bijection $\phi:B\rightarrow C$ (removing eventually a finite subset from $B$) such that
$\sum_{p\in B} \lvert p^{-1}-\phi(p)^{-\beta}\rvert <\infty$. Then, \eqref{eq3.2} implies that
\begin{equation*}
\sum\limits_{p\in B} \bigg| \frac{\sin^2 (t\log p)}{p} -\frac{\sin^2 (\beta t\log \phi (p))}{\phi (p)^\beta} \bigg| < \infty
\end{equation*}
for all $t\in\R$, and hence $G_B=G_{C^\beta}$.

Let $y_n=y_n (s,1)$ be the sequence defined by \eqref{eq3.3} for some fixed constant $s>1$.
Denote $B_n=B\cap (y_n,y_{n+1}]$ and $m_n =\lvert B_n\rvert$. Using Corollary 3.5, we have
\begin{equation}\label{eq3.4}
m_n \leqslant \alpha_n (s,1) =\pi(y_{n+1})-\pi(y_n) \sim \frac{y_n \log^s n}{n(\log \log y_n)^s} \sim \frac{y_n}{n} .
\end{equation}

Set $x_n=y_n(s,\beta^{-1})$; thus $x_n^\beta =y_n$ and $\lim_n x_n/y_n =\infty$. By Corollary 3.5,
\begin{equation}\label{eq3.5}
\pi(x_{n+1})-\pi(x_n) \sim \frac{x_n}{\beta \log x_n (\log \log x_n)^s} \sim \frac{x_n}{n} .
\end{equation}

By \eqref{eq3.4} and \eqref{eq3.5} and since $\lim_n x_n/y_n =\infty$, we get
\begin{equation*}
\lim_n (\pi(x_{n+1})-\pi(x_n)) /m_n =\infty .
\end{equation*}
Hence there exists $n_0$ such that $\lvert \PP \cap (x_n,x_{n+1}] \rvert \geqslant m_n$ for all $n\geqslant n_0$. For each
such $n$ we select a subset $C_n$ of $\PP \cap (x_n,x_{n+1}]$ with $m_n =\lvert B_n\rvert =\lvert C_n \rvert$, then set
$C=\bigcup_{n\geqslant n_0} C_n$ and define a bijection $\phi:B\rightarrow C$ such that $\phi(B_n)=C_n$, for $n\geqslant n_0$.
For any $p\in B_n$ and $q\in C_n$, Lemma 3.2 provides
\begin{equation*}
\bigg| \frac{1}{p}-\frac{1}{q^\beta}\bigg| \leqslant \frac{y_{n+1}-y_n}{y_n^2} \sim \frac{1}{y_n \log^s n} .
\end{equation*}
Henceforth, $\sum_{p\in B} \lvert p^{-1}-\phi(p)^{-\beta}\rvert$ is dominated by $\sum_n m_n /y_n \log^s n$ and by
$\sum_n (n\log^s n)^{-1} <\infty$.
\end{proof}

It is tempting to try to prove that $\beta_0$-representable groups are $\beta$-representable for any $\beta<\beta_0<1$,
using $y_n(s,\beta_0)$ instead of $y_n(s,1)$. However, when $\beta_0 <1$ the previous proof only shows that
$\sum_{p\in B} \lvert p^{-\beta_0}-\phi(p)^{-\beta}\rvert$ is dominated by
$\sum_n (n\log^s n)^{-1} y_n(s,1/\beta_0 -1)$, which is divergent. Nevertheless, replacing the sequences
$y_n(s,c)$ by a sequence with polynomial growth and using \eqref{eq2.2} and \eqref{eq2.3}, we can prove the following.

\begin{thm3.7*}
\emph{If $c_0=0.535$ and $\beta_0 \in (0,1)$, then the $\beta_0$-representable groups are
$\beta$-representable for all $0<\beta <\min (\beta_0,\beta_0 (1-c_0)/(1-\beta_0))$.}
\end{thm3.7*}

\begin{proof}
Let $G=G_{B^{\beta_0}}$ for some subset $B$ of $\PP$ and $\beta_0$, $\beta$ as above. Since
$(n+1)^c -n^c \sim cn^{c-1}$ for all $c>0$, \eqref{eq2.2} and \eqref{eq2.3} yield
\begin{equation}\label{eq3.6}
\pi((n+1)^c) -\pi(n^c) \ll n^{c-1}/\log n \quad \mbox{\rm for all $c>1$,}
\end{equation}
\begin{equation}\label{eq3.7}
n^{c-1}/\log n \ll \pi((n+1)^c) -\pi(n^c) \quad \mbox{\rm for all $c>0$ and $1-c^{-1} >c_0$.}
\end{equation}

Since $\beta < \beta_0(1-c_0)/(1-\beta_0)$, we can choose $a\in (\beta/(1-c_0),\beta_0/(1-\beta_0))$ such that
$a>\beta_0$. It follows that $1-\beta/a>c_0$. We set $B_n =B \cap (n^{a/\beta_0},(n+1)^{a/\beta_0}]$, then use
\eqref{eq3.6} and \eqref{eq3.7} to get
\begin{equation*}
\begin{split}
m_n & = \lvert B_n\rvert \leqslant \pi((n+1)^{a/\beta_0})-\pi(n^{a/\beta_0}) \ll n^{a/\beta_0-1}/\log n ,\\
& n^{a/\beta-1} /\log n \ll \pi( (n+1)^{a/\beta}) -\pi (n^{a/\beta}) .
\end{split}
\end{equation*}

Since $a/\beta_0 < a/\beta$, there exists $n_0>0$ with $m_n \leqslant \lvert \PP \cap (n^{a/\beta},(n+1)^{a/\beta}]\rvert$
for all $n\geqslant n_0$. We select for each such $n$ a subset $C_n$ of $\PP \cap (n^{a/\beta},(n+1)^{a/\beta}]$ with
$\lvert C_n\rvert =m_n$, take $C=\bigcup_{n\geqslant n_0} C_n$, $D=\bigcup_{n\geqslant n_0} B_n$, and define a bijection
$\phi:D\rightarrow C$ such that $\varphi(B_n)=C_n$, obtaining
\begin{equation*}
\sum\limits_{p\in D} \bigg| \frac{1}{p^{\beta_0}} -\frac{1}{\varphi(p)^\beta} \bigg| \leqslant\sum\limits_n
m_n \frac{(n+1)^a-n^a}{n^{2a}} \sim \sum\limits_n \frac{m_n}{n^{a+1}} .
\end{equation*}
The last series is dominated by
\begin{equation*}
\sum\limits_n \frac{n^{a/\beta_0-1}}{n^{a+1}\log n} =
\sum\limits_n \frac{n^{a(1/\beta_0-1)-2}}{\log n} ,
\end{equation*}
which is convergent since $a(1/\beta_0-1)-2<-1$. Using \eqref{eq3.2} also, we get
\begin{equation*}
\sum\limits_{p\in D} \bigg| \frac{\sin^2(\beta_0 t\log p)}{p^{\beta_0}} -
\frac{\sin^2 (\beta t\log \phi(p))}{\phi(p)^\beta} \bigg| < \infty
\end{equation*}
for all $t\in \R$; hence $G_{B^{\beta_0}} =G_{D^{\beta_0}}=G_{\phi(D)^\beta} =G_{C^\beta}$.
\end{proof}

The previous theorem can be sharpened by using the following (see \cite[Lemma 2.13]{AW} or \cite{Bu}).

\begin{lem3.8*}
\emph{Let $(M_\nu)_\nu$ be type $\operatorname{I}_\infty$ factors and $\phi_\nu$, $\psi_\nu$ be normal states
on $M_\nu$ with lists of eigenvalues $\lambda_{\nu,0} \geqslant \lambda_{\nu,1} \geqslant \ldots \geqslant 0$
and respectively $\mu_{\nu,0}\geqslant \mu_{\nu,1}\geqslant \ldots \geqslant 0$. Assume that}
\begin{equation*}
\sum\limits_{\nu} \bigg( 1-\sum\limits_k \sqrt{\lambda_{\nu,k}\mu_{\nu,k}} \bigg) <\infty.
\end{equation*}
\emph{Then $\bigotimes_\nu (M_\nu,\phi_\nu) = \bigotimes_\nu (M_\nu,\psi_\nu)$.}
\end{lem3.8*}

We will also use the following lemma.

\begin{lem3.9*}
\emph{For all $x,y\in [0,1)$ such that $x+y <1$, the following holds:}
\begin{equation*}
1-\sqrt{xy} -\sqrt{(1-x)(1-y)} \leqslant (\sqrt{x}-\sqrt{y})^2 +\sqrt{xy} (\sqrt{x}-\sqrt{y})^2.
\end{equation*}
\end{lem3.9*}

\begin{proof}
Set $a=x+y \geqslant 2\sqrt{b}=2\sqrt{xy}$. The inequality which has to be proved is equivalent to
\begin{equation*}
1-\sqrt{b}-\sqrt{1-a-b} \leqslant (a-2\sqrt{b})(1+\sqrt{b})
\end{equation*}
or
\begin{equation*}
1-a+\sqrt{b} \leqslant \sqrt{b}(a-2\sqrt{b}) +\sqrt{1-a+b} .
\end{equation*}
This follows from $1-a+\sqrt{b}\leqslant \sqrt{1-a+b}$, which is implied by
$(1-a+\sqrt{b})^2 \leqslant 1-a+b$ or equivalently $1-a+2\sqrt{b} \leqslant 1$, which is clear.
\end{proof}

We derive as a corollary the following isomorphism criterion for ITPFI factors of type $M_{\beta,\SSS}$.

\begin{prop3.10*}
\emph{Let $\beta,\beta^\prime \in (0,1]$ and $B$ and $C$ be (infinite) subsets of $\PP$.
If there exists a bijection $\phi:B\rightarrow C$ such that}
\begin{equation*}
\sum\limits_{p\in B} \bigg( \frac{1}{\sqrt{p^\beta}} -\frac{1}{\sqrt{\phi(p)^{\beta^\prime}}}\bigg)^2 < \infty ,
\end{equation*}
\emph{then the factors $M_{\beta,B}$ and $M_{\beta^\prime,C}$ are isomorphic.}
\end{prop3.10*}

\begin{proof}
According to Lemma 3.8 we have to show that
\begin{equation*}
s=\sum\limits_{p\in B} \bigg( 1-\sum\limits_{k\geqslant 0}
\sqrt{p^{-k\beta} (1-p^{-\beta} ) \phi(p)^{-k\beta^\prime} (1-\phi(p)^{-\beta^\prime})}\bigg) < \infty.
\end{equation*}

A direct computation yields
\begin{equation*}
\begin{split}
s & = \sum\limits_{p\in B} \bigg( 1-\frac{\sqrt{(1-p^{-\beta})(1-\phi(p)^{-\beta^\prime})}}{1-\sqrt{p^{-\beta}(\phi(p)^{-\beta^\prime}}} \bigg) \\ &
\sim \sum\limits_{p\in B} \bigg( 1-\sqrt{p^{-\beta}\phi(p)^{-\beta^\prime}}
-\sqrt{(1-p^{-\beta})(1-\phi(p)^{-\beta^\prime})} \bigg) ,
\end{split}
\end{equation*}
and the conclusion follows from the previous lemma.
\end{proof}

Theorem 3.7 can now be improved as follows.

\begin{thm3.11*}
\emph{Let $\beta_0 \in (0,1)$ and $0<\beta < \min (\beta_0,2\beta_0(1-c_0)/(1-\beta_0))$, where $c_0=0.535$.
Then all $\beta_0$-representable groups are $\beta$-representable. In particular, if
$\beta_0 \in (0.07,1)$, then the $\beta_0$-representable groups are $\beta$-representable for any $0<\beta<\beta_0$.}
\end{thm3.11*}

The proof proceeds as in Theorem 3.7. The only difference is that we check that if $a\in (\beta/(1-c_0),2\beta_0/(1-\beta_0))$
and $a>2\beta_0$, then $\sum_{p\in D} (p^{-\beta_0/2}-\phi(p)^{-\beta/2})^2$ is dominated by
\begin{equation*}
\sum\limits_n m_n \frac{((n+1)^{a/2}-n^{a/2})^2}{n^{2a}} \sim \sum\limits_n \frac{m_n}{n^{a+2}} ,
\end{equation*}
and hence by
\begin{equation*}
\sum\limits_n \frac{n^{a(\beta_0^{-1}-1)-3}}{\log n} < \infty .
\end{equation*}

By the previous proposition $M_{\beta_0,B}$ and $M_{\beta,C}$ are isomorphic and in particular one gets
$\operatorname{T}(M_{\beta_0,B})=\operatorname{T}(M_{\beta,C})$.

If Riemann's hypothesis (or something weaker, for example, Lindel\" of's hypothesis or the density hypothesis)
were true, then we could replace $0,07$ by $0$ in Theorem 3.11.

\section{Representable groups and homotheties}
The action by homotheties of $\R^*$ on $\R$ induces a natural action of $\R^*$ on the subgroups of $\R$.
In this section we prove that the class of representable groups is invariant under this action. We also study
the behaviour of $\beta$-representable groups under homotheties.

Denote by $c_\infty$ the set of all sequences $a=\{ a_n\}_n$ with $\lim_n a_n =\infty$ and
$a_n \leqslant a_{n+1}$ for all $n$. Call the elements of $c_\infty$ \emph{divisors}.
Considering the elements of $a$ which repeat, we write $a=a^{\prime m}$, where $a$ and $a^\prime$ coincide as sets,
but $a^\prime$ is a strictly increasing sequence and $m_n \in \N^*$ stands for the number of times each term of $a$ repeats.
The positive integers $m_n$ are called the \emph{multiplicities} of $a$. In other words, if $a$ is the sequence
\begin{equation*}
\underbrace{a_1,\ldots,a_1}_{m_1\, \operatorname{times}} , \underbrace{a_2,\ldots,a_2}_{m_2\, \operatorname{times}  }
,\ldots , \quad \mbox{\rm with $a_1 < a_2 < \ldots ,$}
\end{equation*}
then $a^\prime$ will be the sequence $a_1,a_2,\ldots$, where each $a_k$ appears only once.

For each $a\in c_\infty$ and $t\in\R$, define, as in the previous section,
\begin{equation*}
f_a (2t)=\sum\limits_n \frac{\sin^2 (t\log a_n)}{a_n} =\sum\limits_n m_n \frac{\sin^2 (t\log a_n^\prime)}{a_n^\prime}
\end{equation*}
and consider the subgroup $G(a)=\{ t\in \R: f_a (t)< \infty\}$ of $\R$.
Two sequences $a$ and $b$ from $c_\infty$ are said to be equivalent (we write $a\sim b$) if there exists a bijection
$ \phi:\N^*\rightarrow \N^*$ such that $\sum_n \lvert a_n^{-1} -b_{\phi(n)}^{-1} \rvert < \infty$. It is clear that
$\sim$ is an equivalence relation on $c_\infty$.  The equivalence class of a divisor $a$ is denoted by $[a]$.
The divisors which are equivalent with strictly increasing sequences of primes are called \emph{principal}.
For each $\beta >0$, the divisors equivalent to subsets of $\PP^\beta$, that is, of the form
$\{ p_n^\beta\}_n$, with $p_n < p_{n+1}$, $p_n \in \PP$, are called \emph{$\beta$-principal}.

If $a\sim b$, then $\sum_n a_n^{-1} \sim \sum_n b_n^{-1}$ and according to \eqref{eq3.2}, $f_a(t)\sim f_b(t)$ for
all $t\in \R$; hence $G(a)=G(b)$. We denote $G([a])=G(a)$.

Set ${\mathscr V}=\{ [a]: a\in c_\infty, \sum_n a_n^{-1} =\infty\}$. According to the definition of $\sim$, all sequences
$a=\{ a_n\}_n$ from $c_\infty$ with $\sum_n a_n^{-1} <\infty$ are equivalent and we denote their equivalence class
by ${\mathbf 0}$ and take ${\mathscr V}_0 ={\mathscr V} \cup \{ {\mathbf 0}\}$. We have $G({\mathbf 0})=\R$.
According to the remark at the end of \S 2, the groups $G([a])$, with $[a]\in {\mathscr V}$, are $K_\sigma$ subsets of $\R$.

For $c>0$, $s>1$, we consider the sequence $y(s,c)=\{ y_n (s,c)\}_n$ defined by \eqref{eq3.3}.

\begin{lem4.1*}
\emph{Let $s>1$ and $x=\{ x_n\}_n$ be a sequence such that $x_n < x_{n+1}$, $\lim_n x_n=\infty$,
$\sum_n x_n^{-1} <\infty$ and there exists $c>0$ such that}
\begin{equation}\label{eq4.1}
\frac{x_{n+1}-x_n}{x_n} < \frac{c}{(\log \log x_n)^s} \quad \mbox{\em for all $n\in \N^*$.}
\end{equation}
\emph{Then, for any principal divisor $d$, there exists $m=(m_n)_n$ with $m_n\in \N^*$ and
$m_n \leqslant \alpha_n =\pi(x_{n+1}) -\pi(x_n)$, such that $d \sim x^m$.}
\end{lem4.1*}

\begin{proof}
Let $d\sim B=\{ p_1<p_2 <\ldots\} \subset \PP$. Since $\sum_n x_n^{-1} <\infty$, we may eventually replace
$B$ by $B \cup \{ q_n\}_n$, where $q_n \in \PP$ and $x_n < q_n \leqslant x_{n+1}$ if
$B \cap (x_n,x_{n+1}] =\emptyset$. So we may assume $m_n =\lvert B \cap (x_n,x_{n+1}]\rvert \geqslant 1$.

For any $p\in B$, let $x_n=x_{n(p)}$ be the largest term of $x$ which does not exceed $p$.
Since $\log\log x_n \leqslant \log\log p \leqslant \log\log x_{n+1}$ and $x_{n+1} < 2x_n$ for large $n$, we get
$\lim_n \log\log p/\log \log x_n=1$. On the other hand,
\begin{equation*}
p-x_n < x_{n+1}-x_n < \frac{cx_n}{(\log\log x_n)^s};
\end{equation*}
therefore there exists a constant $c_1 >c$ such that for all $p\in B$,
\begin{equation*}
0\leqslant \frac{1}{x_n} -\frac{1}{p} =\frac{p-x_n}{px_n} <
\frac{c}{p(\log\log x_n)^s} \leqslant \frac{c_1}{p(\log \log p)^s} .
\end{equation*}

Consequently,
\begin{equation*}
\sum\limits_{p\in B} \bigg| \frac{1}{p}-\frac{1}{x_{n(p)}} \bigg| \leqslant c_1
\sum\limits_{p\in B} \frac{1}{p(\log\log p)^s} \leqslant
c_1 \sum\limits_{p\in \PP} \frac{1}{p(\log\log p)^s} < \infty .
\end{equation*}

Taking $\phi (p)=n(p)$, with $p\in B$, we conclude that $d\sim \{ x_n^{m_n}\}_n$ and
$1\leqslant m_n =\lvert B \cap (x_n,x_{n+1}] \rvert \leqslant \alpha_n$.
\end{proof}

\begin{rem*}
By Lemma 3.2, the sequence $\{ y_n(s,c)\}_n$ defined by \eqref{eq3.3} satisfies \eqref{eq4.1} for any
$s>1$ and $c>0$.
\end{rem*}

Let $(\lambda_n)_{n\geqslant 1}$ be a sequence in $(0,1)$. Consider the states $\phi_n =\operatorname{Tr} (\,\cdot\, h_n)$
on $M_2(\C)$, where $h_n$ is the diagonal $2\times 2$ matrix with entries $1/(\lambda_n +1)$ and $\lambda_n/(\lambda_n +1)$.
Let $L_n\in \N^*$ and $M(L_n,\lambda_n)$ be the ITPFI$_2$ factor $\bigotimes_{n\geqslant 1} (M_2(\C),\phi_n)^{\otimes L_n}$.
A computation similar to the one in Lemma 2.6 shows that if $\lim_n \lambda_n =\infty$, then
\begin{equation*}
\operatorname{T}(M(L_n,\lambda_n)) =\bigg\{ t\in \R: \sum\limits_{n} \frac{L_n}{\lambda_n}
\sin^2 \bigg( \frac{t\log\lambda_n}{2}\bigg) < \infty \bigg\} .
\end{equation*}

By Lemma 4.1, $G_B=\{ 2t\in \R: \sum_n m_n x_n^{-1} \sin^2 (t\log x_n)< \infty \}$; hence we obtain the following corollary
which shows that $\operatorname{T}(M_B)$ coincides with the T-group of an ITPFI$_2$ factor whose eigenvalues decrease to zero.

\begin{cor4.2*}
\emph{If $\{ x_n\}_n$ is a sequence as in Lemma {\rm 4.1} and $B\subset \PP$, $B_n=B \cap (x_n,x_{n+1}]$, $m_n=\lvert B_n \rvert$, then}
\begin{equation*}
\operatorname{T}(M_B)=\operatorname{T}(M(m_n,x_n^{-1})).
\end{equation*}
\end{cor4.2*}

This suggests that the factors $M_{1,B}$ may be isomorphic to ITPFI$_2$ factors. Actually it is not hard to see that
this is indeed the case.

\begin{thm4.3*}
\emph{If $\{ x_n\}_n$ is a sequence as in Lemma {\rm 4.1} and $B\subset \PP$, then}
\begin{equation*}
M_{1,B} =M(m_n,x_n^{-1}) .
\end{equation*}
\end{thm4.3*}

\begin{proof}
We have $M(m_n,x_n^{-1}) =\bigotimes_{n\geqslant 1} (\BB (\ell^2 (\N)),\psi_n)^{\otimes m_n}$, where $\psi_n$ has
eigenvalues $\mu_{n,k}=x_n/(x_n+1)$ if $k=0$, $\mu_{n,k}=1/(x_n +1)$ if $k=1$ and $\mu_{n,k}=0$ if $k\geqslant 2$.
We shall estimate $s_B=\sum_{p\in B} (1-\sum_{k\geqslant 0} \sqrt{\lambda_{p,k} \mu_{\phi(p),k}})$ and prove its convergence.

Lemma 3.8 will then imply $M_{1,B}=M(m_n,x_n^{-1})$. By a direct computation
\begin{equation*}
\begin{split}
s_B & = \sum\limits_n \sum\limits_{p\in B_n}
\Bigg( 1-\sqrt{\frac{x_n}{x_n+1} \bigg( 1-\frac{1}{p}\bigg)} -
\sqrt{\frac{1}{p(x_n+1)} \bigg( 1-\frac{1}{p}\bigg)} \Bigg) \\
& = \sum\limits_n \sum\limits_{p\in B_n}
\Bigg( 1-\sqrt{\bigg( 1-\frac{1}{p}\bigg)\frac{x_n}{x_n+1}}  \bigg( 1+\frac{1}{\sqrt{px_n}}\bigg) \Bigg) \\
& \sim t_B =\sum\limits_n \sum\limits_{p\in B_n}
\Bigg( 1-\bigg( 1-\frac{1}{p}\bigg) \bigg( 1-\frac{1}{x_n+1}\bigg) \bigg( 1+\frac{1}{\sqrt{px_n}}\bigg)^2 \Bigg).
\end{split}
\end{equation*}

Since $\lim_n x_{n+1}/x_n=1$ and $\sum_{p\in\PP} p^{-2}<\infty$ we get further
\begin{equation*}
\begin{split}
t_B \sim r_B & = \sum\limits_n \sum\limits_{p\in B_n} \bigg( \frac{1}{p} +\frac{1}{x_n+1} -\frac{2}{\sqrt{px_n}} \bigg) \\
& = \sum\limits_n \sum\limits_{p\in B_n} \Bigg( \bigg( \frac{1}{\sqrt{p}} -\frac{1}{\sqrt{x_n+1}} \bigg)^2
+\frac{2}{\sqrt{p(x_n+1)}}-\frac{2}{\sqrt{px_n}}\bigg) .
\end{split}
\end{equation*}

But
\begin{equation*}
\begin{split}
\sum\limits_n \sum\limits_{p\in B_n} & \bigg( \frac{1}{\sqrt{px_n}}-\frac{1}{\sqrt{p(x_n+1)}}\bigg)
=\sum\limits_n \sum\limits_{p\in B_n} \frac{1}{\sqrt{p(x_n+1)}} \Bigg( \sqrt{\frac{x_n+1}{x_n}} -1 \Bigg) \\
& \leqslant \sum\limits_n \sum\limits_{p\in B_n} \frac{1}{\sqrt{p(x_n+1)}} \bigg( \frac{x_n+1}{x_n}-1\bigg)
\sim \sum\limits_n \sum\limits_{p\in B_n} \frac{1}{x_n \sqrt{px_n}} \sim \sum\limits_{p\in B} \frac{1}{p^2} < \infty .
\end{split}
\end{equation*}

Then, using also the inequality $(\sqrt{a}-\sqrt{b})^2 \leqslant \lvert a-b\rvert$, we get
\begin{equation*}
\begin{split}
s_B & \sim \sum\limits_n \sum\limits_{p\in B_n} \bigg( \frac{1}{\sqrt{p}} -\frac{1}{\sqrt{x_n+1}}\bigg)^2
\leqslant \sum\limits_n \sum\limits_{p\in B_n} \bigg| \frac{1}{p}-\frac{1}{x_n+1}\bigg| \\
& \leqslant \sum\limits_n \sum\limits_{p\in B_n} \Bigg( \bigg( \frac{1}{p} -\frac{1}{x_n}\bigg) +\frac{1}{x_n(x_n+1)}\Bigg)
\sim \sum\limits_n \sum\limits_{p\in B_n} \bigg( \frac{1}{p}-\frac{1}{x_n}\bigg)
+\sum\limits_{p\in B} \frac{1}{p(p+1)}                               < \infty . \qedhere
\end{split}
\end{equation*}
\end{proof}

\begin{lem4.4*}
\emph{Let $y_n=y_n(s,c)$ as in \eqref{eq3.3}, $m_n\in \N^*$ and $n_0\geqslant 1$ such that}
\begin{equation}\label{eq4.2}
m_n \leqslant \frac{cy_n}{2(\log y_n)(\log\log y_n)^s},\quad \mbox{\em for $n\geqslant n_0$.}
\end{equation}
\emph{Then $y^m$ is a principal divisor.}
\end{lem4.4*}

\begin{proof}
By Corollary 3.5, $m_n < \pi(y_{n+1})-\pi(y_n)$ for large $n$.
We choose a subset $B_n$ of $\PP \cap (y_n,y_{n+1}]$ such that $\lvert B_n\rvert =m_n$, take
$B=\bigcup_n B_n$ and check that $\sum_n \sum_{p\in B_n} \lvert p^{-1}-y_n^{-1}\rvert <\infty$, as
in the proof of Lemma 4.1; hence $\{ p\}_{p\in B} \sim y^m$.
\end{proof}

\begin{lem4.5*}
\emph{Let $r_n,s_n\in \N^*$ such that}
\begin{equation*}
M=\sup\limits_n \frac{\lvert r_n-s_n\rvert (\log y_n)(\log\log y_n)^{s+2}}{y_n} < \infty .
\end{equation*}
\emph{Then the divisors $y^r$ and $y^s$ are equivalent.}
\end{lem4.5*}

\begin{proof}
This follows from
\begin{equation*}
\sum\limits_n \bigg| \frac{r_n-s_n}{y_n} \bigg| \leqslant M \sum\limits_n
\frac{1}{(\log y_n)(\log\log y_n)^{s+2}} \sim \sum\limits_n \frac{1}{n\log^2 n} <\infty. \qedhere
\end{equation*}
\end{proof}

\begin{thm4.6*}
\emph{If $G$ is a representable subgroup of $\R$, then $\lambda G$ is representable for any $\lambda\in\R^*$.}
\end{thm4.6*}

\begin{proof}
Let $G=G_B$ (for some $B\subset \PP$) be a representable subgroup of $\R$. Let $\lambda >0$.
We fix $s>1$ and take $y_n=y_n(s,1)$ as in \eqref{eq3.3}. By Lemma 4.1, there exists a divisor $y^m$ equivalent to
$\{ p\}_{p\in B}$ such that $1\leqslant m_n \leqslant \alpha_n =\pi(y_{n+1})-\pi(y_n)$.
By Corollary 3.5 there exists $n_0\geqslant 1$ such that, for all $n\geqslant n_0$,
\begin{equation*}
k_n=\max \bigg( m_n,\frac{y_n}{(\log y_n)(\log\log y_n)^{s+2}} \bigg) \leqslant \alpha_n.
\end{equation*}

Since $k_n-m_n \leqslant (\log y_n)^{-1} (\log\log y_n)^{-s-2} y_n$, Lemma 4.5 implies that
$y^m \sim y^k$ and we may replace $m_n$ by $k_n$. Thus assume that, for any $n$,
\begin{equation}\label{eq4.3}
\alpha_n \geqslant m_n > \frac{y_n}{(\log y_n)(\log\log y_n)^{s+2}} .
\end{equation}

In particular $\lim_n m_n=\infty$. Take $z_n =e^{n/(\lambda \log^s n)} =y_n (s,\lambda^{-1})$, so $z_n^\lambda=y_n$.

Consider $m_n^\prime =[m_n z_n/y_n]$ and the divisor $z^{m^\prime}$. Since
\begin{equation*}
\frac{m_nz_n}{y_n} > \frac{z_n}{(\log y_n)(\log\log y_n)^{s+2}} ,
\end{equation*}
it follows that $\lim_n m_nz_n/y_n =\infty$ and $\lim_n m_n^\prime y_n/m_nz_n =1$. Thus the ratio of the
(non-zero) terms of the series
\begin{equation*}
f_{z^{m^\prime}} (\lambda t) =\sum\limits_n \frac{m_n^\prime \sin^2 (\lambda t\log z_n)}{z_n}
\end{equation*}
and
\begin{equation*}
f_{y^m}(t) =\sum\limits_n \frac{m_n \sin^2 (t\log y_n)}{y_n}
=\sum\limits_n \frac{m_n \sin^2 (\lambda t\log z_n)}{y_n}
\end{equation*}
tends to $1$, implying that $f_{m^\prime}(\lambda t) \sim f_{y^m}(t)$ for all $t\in\R^*$ and therefore
$\lambda G=\lambda G(y^m)=G(z^{m^\prime})$. By \eqref{eq4.3} and Corollary 3.5, there exist constants
$c_1=c_1(s, \lambda)$ and $c_2=c_2(s,\lambda)>0$ such that
\begin{equation*}
\begin{split}
m_n^\prime & \leqslant \frac{m_nz_n}{y_n} \leqslant \frac{\alpha_nz_n}{y_n} \leqslant
\frac{c_1z_n}{(\log y_n)(\log\log y_n)^s} \\
& = \frac{c_1 z_n}{\lambda (\log z_n)(\log\lambda +\log\log z_n)^s} \leqslant\frac{c_2 z_n}{(\log z_n)(\log\log z_n)^s} .
\end{split}
\end{equation*}

We set $m_n^{\prime\prime} =[m_n^\prime /3\lambda c_2]$. Then $\lim_n m_n^{\prime\prime} =\infty$ and
$m_n^{\prime\prime} \leqslant (2\lambda \log z_n)^{-1} (\log\log z_n)^{-s} z_n$ for large $s$
as a result of the previous inequality. Applying Lemma 4.4
to $c=\lambda^{-1}$, we see that $z^{m^{\prime\prime}}$ is a principal divisor. Since $\lim_n m_n^\prime =\infty$, we also have
$\lim_n m_n^\prime /m_n^{\prime \prime} =3\lambda c_2$, which shows that
\begin{equation*}
f_{z^{m^{\prime\prime}}} (t)= \sum\limits_n \frac{m_n^{\prime\prime} \sin^2 (t\log z_n)}{z_n} \sim
\sum\limits_n \frac{m_n^\prime \sin^2 (t\log z_n)}{z_n} =f_{z^{m^\prime}} (t).
\end{equation*}
Thus $\lambda G=G(z^{m^\prime})=G(z^{m^{\prime\prime}})$. Since $z^{m^{\prime\prime}}$ is principal,
we conclude that $\lambda G$ is representable.
\end{proof}

In the remainder of this section we take $y_n(c)=n^c$ for some $c>0$ and consider
the sequence $y(c)=\{ y_n(c)\}_n \in c_\infty$.

\begin{lem4.7*}
\emph{Let $\beta \in (\frac{1}{2},1)$, let $B\subset \PP$ and let $a>2$ such that $\beta a >a-1$.
Then, there exists $m=\{ m_n\}_n$ with $m_n\in\N^*$ and $m_n \leqslant \pi((n+1)^a)-\pi(n^a)$, such that divisors
$\{ p^\beta\}_{p\in B}$ and $y(\beta a)^m$ are equivalent.}
\end{lem4.7*}

\begin{proof}
Set $B_n=B \cap (n^a,(n+1)^a]$. Since $\sum_n n^{-\beta a} <\infty$, we may assume that
$m_n=\lvert B_n \rvert \geqslant 1$. For each $p\in \PP$, let $n_p \in \N^*$ such that
$n_p^a < p\leqslant (n_p+1)^a$ and let $\phi(p)=n_p$, for $p\in B$. Since
$\lim_p n_p^{-1} p^{1/a} =1$ and $(n+1)^c -n^c \sim cn^{c-1}$ for all $c>0$, we find
constants $c_1,c_2 >0$ such that
\begin{equation*}
\bigg| \frac{1}{p^\beta} -\frac{1}{n_p^{\beta a}}\bigg| \leqslant
\frac{(n_p+1)^{\beta a} -n_p^{\beta a}}{p^\beta n_p^{\beta a}} \leqslant
\frac{c_1}{n_p p^\beta} \leqslant \frac{c_2}{p^{\beta +1/a}}
\end{equation*}
for all $p\in\PP$. As a result we get
\begin{equation*}
\sum\limits_{p\in\PP} \bigg| \frac{1}{p^\beta} -\frac{1}{n_p^{\beta a}} \bigg| \leqslant
c_2 \sum\limits_{p\in \PP} \frac{1}{p^{\beta+1/a}} < \infty ;
\end{equation*}
thus $\{ p^\beta\}_{p\in B} \sim y(\beta a)^m$.
\end{proof}

\begin{cor4.8*}
\emph{Let $\beta$, $a$ and $B$ as in Lemma {\rm 4.7}. Then}
\begin{equation*}
\operatorname{T}(M_{\beta,B})=\operatorname{T}(M(L_n,n^{-\beta a})),
\end{equation*}
\emph{where $L_n=\lvert B \cap (n^a,(n+1)^a]\rvert$.}
\end{cor4.8*}

\begin{thm4.9*}
\emph{Let $\beta,B,L_n$ and $a$ be as in Corollary {\rm 4.8}. Then}
\begin{equation*}
M_{\beta,B}=M(L_n,n^{-\beta a}).
\end{equation*}
\end{thm4.9*}

\begin{proof}
We have $M(L_n,n^{-\beta a})=\bigotimes_{n\geqslant 1} (\BB (\ell^2 (\N)),\psi_n)^{\otimes L_n}$, with
the notation from Lemma 4.7 and Corollary 4.8. We also consider the states $\psi_n$ defined by the eigenvalues
$\nu_{n,k} =n^{\beta a}/(n^{\beta a}+1)$ if $k=0$, $\nu_{n,k}=1/(n^{\beta a}+1)$ if $k=1$ and
$\nu_{n,k}=0$ if $k\geqslant 2$. Let $B_n=B \cap (n^a,(n+1)^a]$. According to Lemma 3.8, it is enough to show that
\begin{equation*}
s_B = \sum\limits_n \sum\limits_{p\in B_n} \bigg( 1-\sum\limits_{k\geqslant 0} \sqrt{\lambda_{\beta,p,k}\nu_{n,k}}\bigg) < \infty .
\end{equation*}

We have
\begin{equation*}
\begin{split}
s_B & = \sum\limits_n \sum\limits_{p\in B_n} \Bigg( 1-\bigg( 1+\frac{1}{\sqrt{p^\beta n^{\beta a}}} \bigg)
\sqrt{\bigg( 1-\frac{1}{p^\beta}\bigg) \frac{n^{\beta a}}{n^{\beta a}+1} } \Bigg) \\
& \sim t_B = \sum\limits_n \sum\limits_{p\in B_n} \Bigg( 1-\bigg( 1-\frac{1}{p^\beta}\bigg)
\bigg( 1-\frac{1}{n^{\beta a+1}} \bigg) \bigg( 1+\frac{1}{\sqrt{p^\beta n^{\beta a}}} \bigg)^2 \Bigg) .
\end{split}
\end{equation*}
Since $\lim_n n^{-\beta a} (n+1)^{\beta a} =1$ and $\sum_{p\in\PP} p^{-2\beta} <\infty$, it follows that
\begin{equation*}
t_B=\sum\limits_n \sum\limits_{p\in B_n} \bigg( \frac{1}{p^\beta}
+\frac{1}{n^{\beta a}+1} -\frac{2}{\sqrt{p^\beta n^{\beta a}}} \bigg) =t_B^\prime+t_B^{\prime\prime},
\end{equation*}
where
\begin{equation*}
\begin{split}
t_B^\prime & = 2\sum\limits_n \sum\limits_{p\in B_n} \bigg( \frac{1}{\sqrt{p^\beta n^{\beta a}}}
-\frac{1}{\sqrt{p^\beta (n^{\beta a}+1)}} \bigg) \\ &
\sim \sum\limits_n \sum\limits_{p\in B_n} \frac{1}{p^{\beta/2} n^{\beta a} (\sqrt{n^{\beta a}} +
\sqrt{n^{\beta a}+1})} \sim \sum\limits_n \sum\limits_{p\in B_n} \frac{1}{p^{2\beta}} < \infty
\end{split}
\end{equation*}
and
\begin{equation*}
t_B^{\prime\prime} =\sum\limits_n \sum\limits_{p\in B_n} \bigg( \frac{1}{\sqrt{p^\beta}}
-\frac{1}{\sqrt{n^{\beta a}+1}} \bigg)^2 ,
\end{equation*}
is dominated by
\begin{equation*}
\begin{split}
\sum\limits_n \sum\limits_{p\in B_n} \bigg| \frac{1}{p^\beta} -\frac{1}{n^{\beta a}+1}\bigg| &
\leqslant \sum\limits_n \sum\limits_{p\in B_n} \frac{(n+1)^{\beta a} -n^{\beta a}}{p^\beta n^{\beta a}}
+\sum\limits_n \frac{1}{n^{2\beta a}} \\
& \sim \sum\limits_n \sum\limits_{p\in B_n} \frac{1}{np^\beta} \sim \sum\limits_n \sum\limits_{p\in B_n}
\frac{1}{p^{\beta+1/a}} \leqslant \sum\limits_{p\in\PP} \frac{1}{p^{\beta+1/a}} < \infty .   \qedhere
\end{split}
\end{equation*}
\end{proof}

\begin{lem4.10*}
\emph{Let $a,\beta >0$ be such that $c_0=0.535< 1-a^{-1} <\beta$ and let $m_n\in\N^*$,
$m=\{ m_n\}_n$ be such that $\lim_n n^{1-a} m_n \log n =0$. Then the divisor
$y(\beta a)^m$ is $\beta$-principal.}
\end{lem4.10*}

\begin{proof}
Put $\PP_n =\PP \cap (n^a,(n+1)^a]$. By \eqref{eq3.7} there exists $n_0=n_0(a)$ such that
$m_n \leqslant \lvert \PP_n \rvert$ for all $n\geqslant n_0$; henceforth we can select $B_n \subset \PP_n$
with $\lvert B_n\rvert =m_n$. We take $B=\bigcup_{n\geqslant n_0} B_n$ and $\phi(p)=n$ if $p\in B_n$, getting
\begin{equation*}
\begin{split}
\sum\limits_{p\in B}
\bigg| \frac{1}{p^\beta}
-\frac{1}{\phi(p)^{\beta a}} \bigg| & \leqslant \sum\limits_n \sum\limits_{p\in B_n}
\bigg| \frac{1}{p^\beta}-\frac{1}{n^{\beta a}}\bigg| \leqslant \sum\limits_n \sum\limits_{p\in B_n}
\frac{(n+1)^{\beta a} -n^{\beta a}}{p^\beta n^{\beta a}}
\\
& \sim \sum\limits_n
\sum\limits_{p\in B_n} \frac{1}{np^\beta} \leqslant \sum\limits_{p\in\PP}
\frac{1}{p^{\beta+1/a}} < \infty ;
\end{split}
\end{equation*}
hence $y(\beta a)^m \sim \{ p^\beta\}_{p\in B}$.
\end{proof}

\begin{thm4.11*}
\emph{Let $\beta \in (c_0,1)$ and let $G$ be a $\beta$-representable subgroup of $\R$.
Then $\lambda G$ is $\beta$-representable for any
$\lambda \in ((1-\beta)/(1-c_0),\beta/c_0)$.}
\end{thm4.11*}

\begin{proof}
We take $G=G_{B^\beta}$ for some subset $B\subset\PP$ and set $y_n=y_n (\beta a)=n^{\beta a}$,
$z_n=y_n(\beta a/\lambda)$, $y=\{ y_n\}_n$, $z=\{ z_n\}_n$.

\emph{Case 1:} $1<\lambda <\beta/c_0$. Select $a>2$ such that $c_0 <1-a^{-1}<\beta/\lambda <\beta$ and
keep it fixed throughout the proof. We have $\beta a/\lambda > a-1>1$; hence $\beta a-1>\beta a(1-\lambda^{-1})$.
On the other hand, $\beta a-1>a-2>0$, so we may select $c>0$ such that $\beta a-1>c>\beta a(1-\lambda^{-1})$ and keep it fixed during the proof.
We set $k_n=[n^c] \in \N^*$.

By Lemma 4.7 ($a>2$ and $\beta a>a-1$), we get $m_n\in\N^*$, $m_n \leqslant \alpha_n(a)=\lvert \PP \cap (n^a,(n+1)^a]\rvert$
such that $\{ p^\beta\}_{p\in B} \sim y^m =y(\beta a)^m$. Using $c< \beta a-1<a-1$ and \eqref{eq3.7}
we get $n_0=n_0(a,c)\geqslant 1$ such that $k_n \leqslant \alpha_n (a)$ for all $n\geqslant n_0$.
Moreover, $\sum_n k_n n^{-\beta a} \leqslant \sum_n n^{c-\beta a} <\infty$. Therefore we may replace $m_n$ by
$s_n=\max (m_n,k_n)\leqslant \alpha_n(a)$ and get a new divisor $y^s=y(\beta a)^s$ such that $y^s\sim y^m$
and $[n^c] \leqslant s_n \leqslant \alpha_n(a)$ for all $n\geqslant n_0$. Set
$m_n^\prime =[s_n z_n /y_n]$. We have $s_n z_n/y_n \geqslant c_1 n^{c-\beta a(1-1/\lambda)}$ for some constant
$c_1>0$, so $\lim_n m_n^\prime =\infty$ and $\lim_n m_n^\prime y_n/s_n z_n =1$. As in the proof of Theorem 4.6 we conclude that
$f_{y^s}(t)\sim f_{z^{m^\prime}} (\lambda t)$ for all $t\in\R$, so $\lambda G =\lambda G(y^m)=\lambda G(y^s)=G(z^{m^\prime})$.

It remains to prove that $z^{m^\prime} =y(\beta a/\lambda)^{m^\prime}$ is a $\beta$-principal divisor. This follows
from Lemma 4.10, which applies to $\beta/\lambda$ instead of $\beta$ (since $c_0 < 1-a^{-1} <\beta/\lambda$), and (cf. \eqref{eq3.6})
\begin{equation*}
\frac{m_n^\prime \log n}{n^{a-1}} \leqslant \frac{s_n z_n \log n}{n^{a-1}y_n} \leqslant \alpha_n(a)
\frac{z_n \log n}{n^{a-1} y_n} \leqslant c_2 n^{\beta a(1/\lambda -1)}
\end{equation*}
for some constant $c_2 >0$; hence $\lim_n n^{1-a} m_n^\prime \log n=0$.

\emph{Case 2:} $(1-\beta)/(1-c_0)<\lambda <1$. Since $1-\beta <(1-\beta)/\lambda < 1-c_0$, we may select $a$ such that
$(1-\beta)/\lambda < a^{-1} <1-c_0$. Then $a>2$, $c_0 < 1-a^{-1} <\beta$ and $c_0<1-\lambda/a<\beta$.
As in the beginning of the proof of Case 1 we find $m_n\in \N^*$, $m_n \leqslant \alpha_n(a)$ such that $\{ p^\beta\}_{p\in B}
\sim y^m =y(\beta a)^m$, so $G=G(y^m)$. Set $m_n^\prime =[m_n z_n/y_n]$. Since we have
$m_n z_n/y_n \geqslant n^{\beta a(1/\lambda -1)}$, we get $\lim_n m_n^\prime =\infty$ and
$\lim_n m_n^\prime y_n/m_n z_n =1$, which implies that $f_{y^m}(t) \sim f_{z^{m^\prime}} (\lambda t)$ for all $t\in\R$, so
$\lambda G =\lambda G(y^m)=G(z^{m^\prime})$.

The fact that the divisor $z^{m^\prime} =y(\beta a/\lambda)^{m^\prime}$ is $\beta$-principal is a consequence of
Lemma 4.10, which applies to $a/\lambda$ instead of $a$ (because $c_0 < 1-\lambda/a<\beta$), and of

\begin{equation*}
\frac{m_n^\prime}{\alpha_n (a/\lambda)} \leqslant \frac{\alpha_n(a) z_n}{\alpha_n (a/\lambda)y_n} =
\frac{\alpha_n(a)}{\alpha_n(a/\lambda)} n^{\beta a(1/\lambda -1)} \leqslant c_3 n^{a(\beta -1)(1/\lambda-1)}
\end{equation*}
for some constant $c_3>0$, with the exponent of $n$ in the last term being less than $0$.
\end{proof}

\section{A separation result for Connes T-groups}
In this section we will only take $\beta=1$. Let $\lambda>0$ be a fixed number (actually we will really need only $\lambda=1$) and $A$ be a subset of
$\lambda\N^*$. We denote
\begin{equation*}
f_A (t)=\sum\limits_{n\in A} \frac{1}{n} \sin^2 (2\pi t\log n)
\end{equation*}
and consider the subgroup $G_A=\{ t\in\R: f_A(t) <\infty\}$ of $\R$.

\begin{def5.1*}
If $A$ is a subset of $\lambda \N^*$ for some $\lambda >0$, then a monotonically decreasing sequence
$\varepsilon :A \rightarrow (0,1)$ is called \emph{$A$-admissible} if the following conditions hold:
\begin{equation}\label{eq5.1}
\lim\limits_{n\in A} n\varepsilon_n =\infty ,
\end{equation}
\begin{equation}\label{eq5.2}
\lim\limits_{n\in A} \varepsilon_n=0 ,
\end{equation}
\begin{equation}\label{eq5.3}
\sum\limits_{n\in A} \frac{\varepsilon_n^2}{n} < \infty ,
\end{equation}
\begin{equation}\label{eq5.4}
\sum\limits_{n\in A} \frac{\varepsilon_n}{n}=\infty .
\end{equation}
\end{def5.1*}

\begin{lem5.2*}
\emph{For any $t_1,\ldots,t_k \in \R$, $\ell,N_0\in\N^*$, there exist $r_0 \geqslant 1$ and integers
$q_1,\ldots,q_{r_0}$ such that}
\begin{itemize}
\item[(i)]
\emph{$2N_0 \leqslant q_1 < \ldots < q_{r_0}$;}
\item[(ii)]
\emph{$\displaystyle \frac{1}{2\ell^k} < \sum\limits_{1\leqslant i\leqslant r_0} \frac{1}{q_i} < \frac{1}{\ell^k}$,}
\item[(iii)]
\emph{$\| q_i t_s \| \leqslant \sqrt{k}/\ell$, for all $1\leqslant i\leqslant r_0$ and all $1\leqslant s\leqslant k$.}
\end{itemize}
\end{lem5.2*}

\begin{proof}
We can assume without loss of generality that $N_0 > \ell^k$. Set
$n_0=4N_0 \ell^k+1$ and divide $[0,1)^k$ into $\ell^k$ disjoint cubes of volume $\ell^{-k}$. By the pigeonhole principle the subset
$\{ ( \{ qt_1\},\ldots,\{ qt_k\} ): 1\leqslant q\leqslant n_0\}$ of $[0,1)^k$ contains at least
$4N_0+1$ elements which belong to the same cube. Thus we find integers
$1\leqslant m_1 < m_2 < \ldots < m_{4N_0} <n_0$ such that
\begin{equation*}
\| m_i t_s\| \leqslant \sqrt{k} /\ell,\quad \mbox{\rm for all $1\leqslant i\leqslant 4N_0$ and all $1\leqslant s\leqslant k$.}
\end{equation*}

Take $\{ q_1,\ldots,q_r\} =\{ m_{2N_0},\ldots ,m_{4N_0}\}$. We have $q_1=m_{2N_0} \geqslant 2N_0$ and
\begin{equation*}
\sum\limits_{1\leqslant i\leqslant r} \frac{1}{q_i} \geqslant \frac{2N_0+1}{q_r} > \frac{2N_0+1}{n_0} =
\frac{2N_0+1}{4N_0 \ell^k +1} > \frac{1}{2\ell^k} .
\end{equation*}

On the other hand, $1/q_i < 1/2\ell^k$ since $q_i \geqslant 2N_0 > 2\ell^k$, for $1\leqslant i\leqslant s$. Therefore
we can eventually remove some terms from the tail and get $1\leqslant r_0 \leqslant r$ such that
\begin{equation*}
\frac{1}{2\ell^k} < \sum\limits_{1\leqslant i\leqslant r_0} \frac{1}{q_i} < \frac{1}{\ell^k} .\qedhere
\end{equation*}
\end{proof}

\begin{lem5.3*}
\emph{For any $t_1,\ldots,t_k\in\R$, $\ell,N_0\in\N^*$, there exist $r\geqslant 1$ and integers
$q_1,\ldots,q_r$ such that}
\begin{itemize}
\item[(i)]
$2N_0 \leqslant q_1 < \ldots <q_r$,
\item[(ii)]
$\displaystyle \frac{1}{2\ell} < \sum\limits_{1\leqslant i\leqslant r} \frac{1}{q_i} < \frac{1}{\ell}$,
\item[(iii)]
\emph{$\| q_i t_s\| \leqslant \sqrt{k}/\ell$, for all $1\leqslant i\leqslant r$ and all $1\leqslant s\leqslant k$.}
\end{itemize}
\end{lem5.3*}

\begin{proof}
Apply Lemma 5.2 and obtain $r_1\geqslant 1$ and integers $q_{r_1} > \ldots > q_1 \geqslant 2N_0$ such that
\begin{equation*}
\frac{1}{\ell^k} > \sum\limits_{1\leqslant i\leqslant r_1} \frac{1}{q_i} > \frac{1}{2\ell^k} ,
\end{equation*}
\begin{equation*}
\| q_i t_s \| \leqslant \sqrt{k}/\ell,\quad \mbox{\rm for all $1\leqslant i\leqslant r_1$ and all $1\leqslant s\leqslant k$.}
\end{equation*}

Taking $N_1>\max (\ell^k,\frac{1}{2} q_{r_1})$ and applying Lemma 5.2 again we find
$r_2>r_1$ and some new integers $q_{r_2} > \ldots > q_{r_1+1} \geqslant 2N_1 > q_{r_1}$ such that
\begin{equation*}
\frac{1}{\ell^k} > \sum\limits_{r_1 < i\leqslant r_{i_2}} \frac{1}{q_i} > \frac{1}{2\ell^k} ,
\end{equation*}
\begin{equation*}
\| q_i t_s \| \leqslant \sqrt{k}/\ell,\quad \mbox{\rm for all $r_1< i\leqslant r_1$ and all $1\leqslant s\leqslant k$.}
\end{equation*}

Applying Lemma 5.2 $\ell^{k-1}$ times with $N_j >\max (\ell^k,\frac{1}{2} q_{r_j})$ at each step $j$, we end up with
$r=r_{\ell^{k-1}}$ and integers $q_r > \ldots > q_1 \geqslant 2N_0$ such that
\begin{equation*}
\frac{1}{\ell} > \sum\limits_{1\leqslant i\leqslant r} \frac{1}{q_i} > \frac{1}{2\ell},
\end{equation*}
\begin{equation*}
\| q_i t_s \| \leqslant \sqrt{k}/\ell,\quad \mbox{\rm for all $1\leqslant i\leqslant r$ and all $1\leqslant s\leqslant k$.} \qedhere
\end{equation*}
\end{proof}

\begin{cor5.4*}
\emph{For any $t_1,\ldots,t_k\in\R$, there exist a subset $A$ of $\N^*$ and an $A$-admissible sequence
$\varepsilon=\{ \varepsilon_n\}_{n\in A}$ such that}
\begin{equation}\label{eq5.5}
\| nt_s\| \leqslant \varepsilon_n \sqrt{k}, \quad \mbox{\em for all $n\in A$ and all $1\leqslant s\leqslant k$.}
\end{equation}
\end{cor5.4*}

\begin{proof}
Let $\{ c_\ell\}_\ell$ be a sequence which decreases monotonically to $0$ and $\ell c_\ell \geqslant 1$, $\lim_\ell \ell c_\ell =\infty$ and
$\sum_\ell c_\ell^2 /\ell <\infty$. Applying Lemma 5.3 for each $\ell \geqslant 1$ we find a strictly increasing sequence of
integers $\{ q_n\}_n$ such that
\begin{equation}\label{eq5.6}
\frac{1}{2\ell} < \sum\limits_{r_\ell < i\leqslant r_{\ell+1}} \frac{1}{q_i} < \frac{1}{\ell} ,
\end{equation}
\begin{equation}\label{eq5.7}
\| q_i t_s \| \leqslant \sqrt{k}/\ell, \quad \mbox{\rm for all $r_\ell <i\leqslant r_{\ell+1}$ and $1\leqslant s\leqslant k$.}
\end{equation}

We take $A=\{ q_n\}_n$ and $\varepsilon_{q_i} =c_\ell$ if $r_\ell < i\leqslant r_{\ell+1}$. Then, condition \eqref{eq5.5}
is clearly fulfilled since
\begin{equation*}
\| q_i t_s\| \leqslant \sqrt{k}/\ell \leqslant c_\ell \sqrt{k} =\varepsilon_{q_i} \sqrt{k},
\quad \mbox{\rm for $r_\ell < i\leqslant r_{\ell+1}$, $1\leqslant s\leqslant k$.}
\end{equation*}

Moreover, \eqref{eq5.1} holds because $q_i c_\ell \geqslant q_{r_\ell} c_\ell \geqslant r_\ell c_\ell \geqslant \ell c_\ell$ if
$i\geqslant r_\ell$, \eqref{eq5.2} is obvious, while \eqref{eq5.3} and \eqref{eq5.4} follow from
\begin{equation*}
\sum\limits_{n\in A} \frac{\varepsilon_n}{n} =\sum\limits_{\ell} \sum\limits_{r_\ell < i\leqslant r_{\ell+1}}
\frac{c_\ell}{q_i} \geqslant \sum\limits_\ell \frac{c_\ell}{2\ell} =\infty\end{equation*}
and
\begin{equation*}
\sum\limits_{n\in A} \frac{\varepsilon_n^2}{n} =\sum\limits_\ell \sum\limits_{r_\ell < i\leqslant r_{\ell+1}}
\frac{c_\ell^2}{q_i} \leqslant \sum\limits_\ell \frac{c_\ell^2}{\ell} < \infty ,
\end{equation*}
respectively.
\end{proof}

The following definition is inspired by Corollary 2.2.

\begin{def5.5*}
Let $\lambda >0$, $A$ be a subset of $\lambda \N^*$, and $\varepsilon$ be an $A$-admissible sequence. A subset
$B\subset \PP$ is called an $\varepsilon$-lifting of $A$ if there exists a surjective map $h:B \rightarrow A$ such that
\begin{equation}\label{eq5.8}
0 < \inf\limits_{n\in A} \frac{n}{\varepsilon_n} \sum\limits_{h(p)=n} \frac{1}{p} \leqslant\sup\limits_{n\in A}
\frac{n}{\varepsilon_n} \sum\limits_{h(p)=n} \frac{1}{p}<\infty ,
\end{equation}
\begin{equation}\label{eq5.9}
\sup\limits_{p\in B} \frac{\lvert h(p)-\log p\rvert}{\varepsilon_{h(p)}} < \infty .
\end{equation}
\end{def5.5*}

\begin{rems*}
(i) If $\varepsilon$ is $A$-admissible and $B\subset \PP$ and $h:B\rightarrow A$ are as above, then
$\sum_{p\in B} p^{-1}=\infty$.

(ii) If \eqref{eq5.2} and \eqref{eq5.9} hold and $A\subset \N^*$, then $h(p)$ is the nearest integer from $\log p$ for large $p$.
\end{rems*}

For any $\lambda >0$, any subset $A\subset \lambda \N^*$, and any $A$-admissible sequence $\varepsilon$ and any
$a\in [0,1)=\R /\Z$, we define
\begin{equation*}
H_a (\varepsilon) =\bigg\{ t\in\R : \sup\limits_{n\in A} \frac{\| nt-a\|}{\varepsilon_n} < \infty \bigg\} .
\end{equation*}

We have $H_{-a}(\varepsilon)=-H_a(\varepsilon)$ and $H_{a_1}(\varepsilon)+H_{a_2}(\varepsilon)
\subset H_{a_1+a_2} (\varepsilon)$. In particular, if $\Gamma$ is a subgroup of $\R$, then
$\bigcup_{a\in\Gamma} H_a(\varepsilon)$ is also a subgroup of $\R$. The features of the sets
$H_a(\varepsilon)$ from the next three lemmas will play an important role in the proof of the separation theorem.

\begin{lem5.6*}
\emph{Let $\lambda >0$, $A$ be a subset of $\lambda\N^*$, $\varepsilon$ be an $A$-admissible sequence,
$B\subset \PP$ an $\varepsilon$-lifting of $A$ and $a\notin \{ 0,\frac{1}{2}\}$ in $\R /\Z$. Then}
\begin{equation*}
H_a(\varepsilon) \cap G_B =\emptyset .
\end{equation*}
\end{lem5.6*}

\begin{proof}
Let $h:B\rightarrow A$ be a surjection which implements the $\varepsilon$-lifting $B\subset \PP$ of $A$ and
$t\in H_a (\varepsilon) \cap G_B$. Let $\delta >0$ be such that the interval
$[\| a\|-2\delta,\| a\|+2\delta ]$ does not intersect $\frac{1}{2}\Z$ and set
\begin{equation*}
\mu=\inf \{ \sin^2 (2\pi x): \| a\|-2\delta \leqslant x\leqslant \| a\|+2\delta \} > 0.
\end{equation*}
By \eqref{eq5.9}, there exists $c_0 >0$ such that for all $p\in B$ with $h(p)=n$ and $n\in A$ we have
\begin{equation*}
\| nt\|+c_0 \varepsilon_n \lvert t\rvert > \| t\log p\| > \| nt\| -c_0\varepsilon_n \lvert t\rvert .
\end{equation*}

Since $t\in H_a(\varepsilon)$, there exists $c_1>0$ such that $\| nt-a\| \leqslant c_1 \varepsilon_n$ for
all $n\in A$. Hence $\| a\|+c_1\varepsilon_n > \| nt\| >\| a\|-c_1 \varepsilon_n$ for all $n\in A$.
But $\lim_{n\in A} \varepsilon_n=0$; thus
there exists $n_1 >0$ such that $\max (c_1 \varepsilon_n,c_0 \varepsilon_n \lvert t\rvert) < \delta$ for
all $n\geqslant n_t$ and we get
\begin{equation*}
\| a\|-2\delta < \| t\log p\| < \| a\| +2\delta,\quad \mbox{\rm for all $n\in A$ and $n\geqslant t$},
\end{equation*}
whence
\begin{equation*}
f_B(t) =\sum\limits_{n\in A} \sum\limits_{h(p)=n} \frac{\sin^2 (2\pi t\log p)}{p} \geqslant \mu
\sum\limits_{n\in A} \sum\limits_{h(p)=n} \frac{1}{p} =\mu\sum\limits_{p\in B} \frac{1}{p} =\infty .
\end{equation*}

This contradicts the assumption $t\in G_B$.
\end{proof}

\begin{lem5.7*}
\emph{If $B$ is an $\varepsilon$-lifting of $A$, then}
\begin{equation*}
H_0(\varepsilon) \cup H_{1/2} (\varepsilon) \subset G_B .
\end{equation*}
\end{lem5.7*}

\begin{proof}
We have $\lvert h(p) -\log p\rvert < c_1 \varepsilon_{h(p)}$ for $p\in B$, for some constant $c_1>0$.
We pick $t\in H_0(\varepsilon)$, evaluate $f_B(t)$ using $\| a+b\|^2 \leqslant 2(\| a\|^2+\| b\|^2)$
for $a,b\in\R$, and get
\begin{equation*}
f_B(t) = \sum\limits_{n\in A} \sum\limits_{h(p)=n} \frac{\sin^2 (2\pi t\log p)}{p} \leqslant
4\pi^2 \sum\limits_{n\in A} \sum\limits_{h(p)=n} \frac{\| t\log p\|^2}{p}
\leqslant 8\pi^2 \sum\limits_{n\in A} \sum\limits_{h(p)=n} \frac{\| tn\|^2+c_1^2t^2\varepsilon_n^2}{p} .
\end{equation*}

Since $t\in H_0(\varepsilon)$, there exists $c_2=c_2(t)>0$ such that $\| nt\| \leqslant c_2\varepsilon_n$
for all $n\in A$. Therefore there exists a constant $c_3=c_3(t)>0$ such that
\begin{equation*}
f_B(t) \leqslant c_3 \sum\limits_{n\in A} \varepsilon_n^2 \sum\limits_{h(p)=n} \frac{1}{p} .
\end{equation*}
But $\sup_{n\in A} (n/\varepsilon_n) \sum_{h(p)=n} p^{-1} < \infty$; thus there exists $c_4>0$
such that $\sum_{h(p)=n} p^{-1} \leqslant c_4 \varepsilon_n/n$ for all $n\in A$ and consequently
\begin{equation*}
f_B(t) \leqslant c_3 c_4 \sum\limits_{n\in A} \frac{\varepsilon_n^3}{n}
\leqslant c_3c_4 \sum\limits_{n\in A} \frac{\varepsilon_n^2}{n} < \infty,
\end{equation*}
which shows that $t\in G_B$.

The proof of $H_{1/2} (\varepsilon) \subset G_B$ is similar, using \eqref{eq5.9}, \eqref{eq3.1} and
\begin{equation*}
\sum\limits_{n\in A} \sum\limits_{h(p)=n} \frac{\sin^2 (2\pi nt)}{p}
\leqslant 4\pi^2 \sum\limits_{n\in A} \sum\limits_{h(p)=n}
\frac{\| nt+\frac{1}{2}\|^2}{p} \leqslant c_6 \sum\limits_{n\in A} \frac{\varepsilon_n^3}{n} < \infty,
\end{equation*}
for some constant $c_6=c_6(t)>0$.
\end{proof}

\begin{lem5.8*}
\emph{Let $\lambda >0$, $A$ be a subset of $\lambda\N^*$ and $\varepsilon$ be an $A$-admissible
sequence. Then, there exists a subset $B\subset \PP$ such that}
\begin{equation*}
H_0(\varepsilon) \cup H_{1/2} (\varepsilon) \subset G_B .
\end{equation*}
\end{lem5.8*}

\begin{proof}
The sets $B_n=\{ p\in\PP: \lvert n-\log p\rvert <\varepsilon_n\}$, with $n\in A$, are mutually
disjoint if $2\varepsilon_n <1$. By Lemma 2.1 and the proof of Corollary 2.2 we get
\begin{equation*}
\lim\limits_{n\in A} \frac{n}{\varepsilon_n} \sum\limits_{p\in B_n} \frac{1}{p} =2.
\end{equation*}
Thus taking $B=\bigcup_{n\in A} B_n$ and $h(p)=n$ if $p\in B_n$, it follows that $B$ is an $\varepsilon$-lifting
of $A$. By Lemma 5.7 we gather that $H_0(\varepsilon) \cup H_{1/2} (\varepsilon) \subset G_B$,
which completes the proof.
\end{proof}

The proof of the following technical statement is a refinement of Corollary 5.4.
In addition we shall perturb at each step the `platoon' $q_{r_\ell +1},\ldots ,q_{r_{\ell+1}}$ to
keep $\| q_i u\|$ far from $0$ and $\frac{1}{2}$.

\begin{prop5.9*}
\emph{Let $t_1,\ldots,t_k \in\R$ be linearly independent over $\Q$,
$G=\Z t_1 +\ldots +\Z t_k$ and $u\in \R \setminus G$. Then, there exists a subset
$B\subset \PP$ such that $G\subset G_B$ and $u\notin G_B$.}
\end{prop5.9*}

\begin{proof}
By Theorem 4.6 we may assume, without loss of generality, that $t_k=1$. Write $u=u^\prime+u^{\prime\prime}$, with $1,t_1,\ldots,t_{k-1},u^{\prime\prime}$
linearly independent over $\Q$ and $u^\prime =\sum_{j=1}^k r_j t_j$, with $r_j=a_j/b_j$, $a_j\in\Z$, $b_j\in\N^*$,
$\gcd (a_j,b_j)=1$ (if $a_j=0$ we take $b_j=1$). Denote
$C=2+\sqrt{k}+\sum_{j=1}^{k-1} (\lvert a_j\rvert +b_j)$. Consider also an $\N^*$-admissible sequence $\{ c_\ell \}_{\ell}$
with $\ell c_\ell \geqslant 1$, for instance $c_\ell =(\log \ell)^{-1}$.

We will prove the existence of a real number $a$ with $a\notin \{ 0,\frac{1}{2}\} \mod \Z$, of a set
$A\subset \N^*$ and of an $A$-admissible sequence $\varepsilon$ such that $t_j \in H_0(\varepsilon) \cup H_{1/2} (\varepsilon)$, for
$1\leqslant j\leqslant k$, and $u\in H_a(\varepsilon)$. Then, this will enable us to produce an $\varepsilon$-lifting $B\subset \PP$
with Lemma 5.8 and to conclude that $t_j \in G_B$ with Lemma 5.7 and that $u\notin G_B$
with Lemma 5.6.

The existence of $a$, $A$ and $\varepsilon$ follows once we succeed in showing the existence of some
$s\in \{ 1,\ldots,k-1\}$ and of a family of finite subsets $\{ Q_\ell\}_{\ell \geqslant 1}$ of $\N^*$ such that
\begin{equation}\label{eq5.10}
\sup Q_\ell < \inf Q_{\ell+1} ,\quad \mbox{\rm for $\ell\in\N^*$,}
\end{equation}
\begin{equation}\label{eq5.11}
\frac{1}{4\ell} < \sum\limits_{q\in Q_\ell} \frac{1}{q} < \frac{1}{\ell} ,
\end{equation}
\begin{equation}\label{eq5.12}
\max \big( \| qu-a\|.\max\limits_{j\leqslant k-1,j\neq s} \| qt_j\|\big) < C \alpha_\ell ,\quad\mbox{\rm for $q\in Q_\ell$, $\ell\in\N^*$,}
\end{equation}
\begin{equation}\label{eq5.13}
\min \big( \| qt_s \| ,\| qt_s -\tfrac{1}{2}\| \big) < C\alpha_\ell ,\quad
\mbox{\rm for $q\in Q_\ell$, $\ell \in \N^*$,}
\end{equation}
then take $A=\bigcup_\ell Q_\ell \subset \N^*$ and $\varepsilon_q =c_\ell$ if $q\in Q_\ell$ (the $A$-admissibility of
$\varepsilon$ follows from \eqref{eq5.11}).

The construction of such $s$ and $\{ Q_\ell\}_{\ell}$  relies on Lemma 5.3 and on Kronecker's theorem. To proceed,
fix $N>0$ and $\ell \in\N^*$.

Consider first the case $u^{\prime\prime} \neq 0$ and take $a=\frac{1}{4}$ and $s$ any element from
$\{ 1,\ldots,k-1\}$. With Kronecker's theorem we find $m_0\in\N^*$ such that
\begin{equation}\label{eq5.14}
\max \bigg( \| m_0 u^{\prime\prime} -\tfrac{1}{4} \|, \max\limits_{j\leqslant k-1} \bigg\| \frac{m_0 t_j}{b_j}\bigg\| \bigg) < c_\ell .
\end{equation}

By Lemma 5.3, there exists a finite set $Q_0 \subset \N^*$ such that
\begin{equation}\label{eq5.15}
\max (2N,m_0) < \inf Q_0 ,
\end{equation}
\begin{equation}\label{eq5.16}
\frac{1}{2\ell} < \sum\limits_{q\in Q_0} < \frac{1}{\ell},
\end{equation}
\begin{equation}\label{eq5.17}
\max ( \| qu\|,\| qt_1\|,\ldots ,\| qt_{k-1}\|) \leqslant c_\ell \sqrt{k},\quad \mbox{\rm for $q\in Q_0$.}
\end{equation}

Take $Q_\ell =\{ q^\prime=q+m_0:q\in Q_0\}$. Then $\frac{1}{2} < q/q^\prime <1$, for $q\in Q_0$, which we combine
with \eqref{eq5.16} to get \eqref{eq5.11}. Moreover, \eqref{eq5.14} and \eqref{eq5.17} provide
\begin{equation*}
\max \Big( \| q^\prime u-a\|,\max\limits_{j\leqslant k-1} \| q^\prime t_j\| \Big) < C \alpha_\ell ,\quad
\mbox{\rm for $q^\prime \in Q_\ell$, $\ell\in \N^*$},
\end{equation*}
so in this case $t_1,\ldots,t_k \in H_0(\varepsilon)$ and $u\in H_{1/4} (\varepsilon)$.

When $u^{\prime\prime}=0$ we choose $s\in \{ 1,\ldots,k-1\}$ such that $r_s \notin \N^*$ and take $a=\frac{1}{2} r_s$.
We find, again with Kronecker's theorem, $m_0\in\N^*$ such that
\begin{equation}\label{eq5.18}
\max \bigg( \bigg\| \frac{m_0 t_s}{b_s} -\frac{1}{2b_s} \bigg\| ,\max\limits_{j\leqslant k-1,j\neq s} \bigg\|
\frac{m_0 t_j}{b_j}\bigg\| \bigg) < c_\ell .
\end{equation}

We then select, using Lemma 5.3, a finite set $Q_0 \subset \N^*$ such that \eqref{eq5.15}, \eqref{eq5.16}, \eqref{eq5.17} are fulfilled and take
$Q_\ell$ as in the first case. Then \eqref{eq5.11} follows for a similar reason, whilst \eqref{eq5.18} and \eqref{eq5.17} provide
\begin{equation*}
\max \Big( \| q^\prime u-\tfrac{1}{2} r_s \|,\| q^\prime t_s -\tfrac{1}{2}\| ,
\max\limits_{j\leqslant k-1,j\neq s} \| q^\prime t_j \| \Big) < C \alpha_\ell, \quad\mbox{\rm for $q^\prime \in Q_\ell$, $\ell \in \N^*$.}
\end{equation*}

In this case $t_s \in H_{1/2} (\varepsilon)$, $t_j \in H_0(\varepsilon)$ if $j\neq s$ and $u\in H_{r_s/2}(\varepsilon)$.
\end{proof}

\begin{cor5.10*}
\emph{Let $t_1,\ldots,t_k\in\R$ be linearly independent over $\Q$, let $G=\Z t_1 +\ldots +\Z t_k$ and let $u_1,\ldots,u_n\in G^c$.
Then, there exists a representable group $\Gamma \subset \R$ such that}
\begin{equation*}
G\subset \Gamma \subset \{ u_1,\ldots,u_n\}^c .
\end{equation*}
\end{cor5.10*}

\begin{proof}
By Proposition 5.9, there exists for each $1\leqslant i\leqslant n$ a representable group $\Gamma_i$ which contains $G$ and does not contain $u_i$.
The intersection of two representable groups is clearly representable since $0\leqslant f_{B_1 \cup B_2} \leqslant f_{B_1} +f_{B_2}$ if
$B_1,B_2 \subset \PP$ and $0\leqslant f_{B_1} \leqslant f_{B_2}$ if $B_1 \subset B_2 \subset \PP$. Thus
$G_{B_1} \cap G_{B_2} =G_{B_1 \cup B_2}$. Therefore $\Gamma =\bigcap_{i=1}^n \Gamma_i$ is a representable group which contains $G$ and does
not contain any $u_i$.
\end{proof}

The following `truncation' lemma allows us to replace `finitely generated' by `countable'.

\begin{lem5.11*}
\emph{Let $G=G_I$ be a representable group with $I$ a subset of $\PP$, let $A\subset G$ and $B\subset G^c$ be finite sets and $N,m,M >0$.
Then, there exists a finite subset $I_0 \subset I$ such that $\inf I_0 \geqslant N$ and}
\begin{equation*}
f_{I_0}(t) <m,\quad \mbox{\em for all $t\in A$,}
\end{equation*}
\begin{equation*}
f_{I_0}(t) >M, \quad \mbox{\em for all $t\in B$.}
\end{equation*}
\end{lem5.11*}

\begin{proof}
Truncate the series twice.
\end{proof}

\begin{thm5.12*}
\emph{Let $G$ be a countable subgroup of $\R$, and $\Sigma$ be a countable
subset of $G^c$. Then, there exists a representable group $\Gamma$ such that}
\begin{equation*}
G \subset \Gamma \subset \Sigma^c.
\end{equation*}
\end{thm5.12*}

\begin{proof}
Let $G=\bigcup_n S_n$, where $S_n$ is an increasing sequence of finite sets. Put
$G_n =\sum_{t\in S_n} \Z t$. Let $\Sigma =\{ h_1,\ldots,h_n,\ldots\}\subset G^c$ and $\Sigma_n =\{ h_1,\ldots, h_n\}$.
With Corollary 5.10 we produce for each $n$ a representable group $\Gamma_n$ such that
$G_n \subset \Gamma_n \subset \Sigma_n^c$. We apply Lemma 5.11 to find a sequence $\{ I_n\}_n$ of finite subsets of $\PP$ such that
$\sup I_n < \inf I_{n+1}$ and
\begin{equation*}
\begin{split}
f_{I_n} (t) & < 1/n^2 ,\quad \mbox{\rm for all $t\in S_n$,} \\
f_{I_n} (t) & > 1 ,\quad \mbox{\rm for all $t\in\Sigma_n$.}
\end{split}
\end{equation*}

Taking $B=\bigcup_n I_n$, we obtain
\begin{equation*}
f_B(t)=\sum\limits_n f_{I_n}(t) \leqslant \sum\limits_{n\leqslant k} f_{I_n}(t) +
\sum\limits_{n>k} \frac{1}{n^2} < \infty
\end{equation*}
for all $t\in S_k$. Hence $f_B(t)< \infty$ for all $t\in \bigcup_k S_k =G$, and therefore
$G\subset \Gamma =G(f_B)$. On the other hand, we have
\begin{equation*}
f_B(t) \geqslant \sum\limits_{n\leqslant n_0} f_{I_n} (t) +\sum\limits_{n>n_0} 1 =\infty
\end{equation*}
for all $t\in \Sigma_k$, showing that $f_B(t)=\infty$ for all
$t\in \bigcup_k \Sigma_k =\Sigma$ and consequently $\Sigma \cap \Gamma =\emptyset$.
\end{proof}

The following are immediate consequences of Theorems 5.12 and 3.6.

\begin{cor5.13*}
\emph{Let $\beta \in (0,1]$, let $G$ be a countable subgroup of $\R$ and let $\Sigma$ be a countable subset of $G^c$.
Then, there exists a subset $S=S(\beta,G,\Sigma) \subset \PP$ such that}
\begin{equation*}
G\subset \operatorname{T}(M_{\beta,S}) \subset \Sigma^c .
\end{equation*}
\end{cor5.13*}

\begin{cor5.14*}
\emph{Let $G$ and $H$ be countable subgroups of $\R$ such that $G\cap H=\{ 0\}$. Then, for any
$\beta \in (0,1]$, there exists a $\beta$-representable group $\Gamma$ such that $G\subset \Gamma$ and
$\Gamma \cap H=\{ 0\}$.}
\end{cor5.14*}

\begin{cor5.15*}
\emph{Let $\beta \in (0,1]$. Then, any countable subgroup of $\R$ is the intersection of all $\beta$-representable
groups which contain it. In particular, all countable subgroups of $\R$ are $\beta$-admissible.}
\end{cor5.15*}

\begin{cor5.16*}
\emph{If $G$ is an admissible subgroup of $\R$, then $\lambda G$ is admissible for any $\lambda\in \R^*$.}
\end{cor5.16*}

A natural problem which arises is to study the automorphisms of a representable group.
In the case of cyclic subgroups of $\R$, known to be representable by Proposition 2.5, the automorphism group is simply $\{ -1,1\}$.
In general. the automorphisms of a representable group which are not homotheties
are complicated and they are not continuous with respect to the topology induced from $\R$.
We will consider in the sequel only homotheties.

If $G$ is representable, $\operatorname{Omo}(G) =\{ \lambda \in\R^*: \lambda G=G\}$ is a subgroup pf
$\operatorname{Aut}(G)$. If $K$ is a number field, denote by $A_K$ its integer ring
and by $U_K$ its group of units.

\begin{prop5.17*}
\emph{For any real number field $K$, there exists a non-zero representable group $\Gamma$ which is an $A_K$-module.}
\end{prop5.17*}

\begin{proof}
Consider an integral basis $w_1,\ldots ,w_r$ of $K$. By Theorem 5.12, there exists a representable
group $G$ which contains $K$.
Since the class of representable groups is closed under finite intersections and homotheties,
$\Gamma =\bigcap_{j=1}^r w_j^{-1} G$ is a representable group which contains $K$. A real number $t$ belongs to $\Gamma$ if and
only if $w_j t\in G$ for all $1\leqslant j\leqslant r$, or equivalently $A_K t \subset G$. Therefore,
if $t\in\Gamma$ and $a\in A_K$, then $A_K at =A_K t \subset G$, which implies $at\in\Gamma$.
\end{proof}

\begin{cor5.18*}
\emph{For any real number field $K$, there exists a non-zero representable group $\Gamma$ such that}
\begin{equation*}
U_K \subset \operatorname{Omo}(\Gamma) .
\end{equation*}
\end{cor5.18*}

\begin{proof}
Let $\Gamma \neq \{ 0\}$ be a representable group which is an $A_K$-module.
If $u\in U_K$, then $u,u^{-1} \in A_K$ and $u\Gamma \subset\Gamma$, $u^{-1} \Gamma \subset \Gamma$,
which show that $u\Gamma =\Gamma$.
\end{proof}

Using the separation part of Theorem 5.12 as well we obtain the following:

\begin{cor5.19*}
(i) \emph{Let $G_0$ be a countable $A_K$-module and $M$ be a countable subset of $G_0^c$.
Then, there exists a representable group $\Gamma$ which is an $A_K$-module and $G_0 \subset \Gamma \subset M^c$.}

(ii) \emph{Let $G_0$ be a countable group such that $U_K G_0 \subset G_0$ and $M$ be a countable subset
of $G_0^c$. Then, there exists a representable group $\Gamma \neq \{ 0\}$ such that
$\Gamma \cap M=\emptyset$ and $U_K \subset \operatorname{Omo}(\Gamma )$.}
\end{cor5.19*}

\begin{cor5.20*}
\emph{Let $P_1,\ldots,P_m$ be monic polynomials in $\Z [X]$ with $P_j(0)\in \{ -1,1\}$ and $b_j$ be real numbers with
$P_j(b_j)=0$ for all $1\leqslant j\leqslant m$. Then, there exists a representable group $\Gamma \neq \{ 0\}$ such that}
\begin{equation*}
b_1,\ldots ,b_m \in \operatorname{Omo} (\Gamma) .
\end{equation*}
\end{cor5.20*}

\begin{proof}
Take $K=\Q (b_1,\ldots ,b_m)$ and notice that $b_1,\ldots ,b_m \in U_K$.
\end{proof}

In order to study admissibility with respect to subgroups of $\R$, we associate to each pair $(E,H)$ of subgroups of $\R$ the abelian group
\begin{equation*}
\DD_{E,H} =\bigg( \bigcap\limits_{G\in \RR_E} (G+H)\bigg) \bigg\slash (E+H),
\end{equation*}
where $\RR_E =\{ G \  \mbox{\rm representable}: E\subset G\}$. The group $\DD_{E,H}$, which we call the \emph{defect
of $E$ modulo $H$}, provides an obstruction for admissibility modulo $H$. The motivation for introducing the defect groups
is that they coincide with some firs cohomology groups. This connection will be explained in a forthcoming paper.
Some properties of defect groups are listed below:
\begin{itemize}
\item[(1)] a group inclusion $H_1 \subset H_2$ induces a group morphism $\DD_{E,H_1} \rightarrow \DD_{E,H_2}$;
\item[(2)] if $E$ is representable, then $\DD_{E,H}=\{ 0\}$ for any $H$;
\item[(3)] $\DD_{E,\{ 0\}}=\{ 0\}$ if and only if $E$ is admissible;
\item[(4)] $\DD_{E,H} =\DD_{E,E+H}$;
\item[(5)] if $E \subset H_i$ and $\DD_{E,H_i}=\{ 0\}$ for all $i$, then $\DD_{E, \cap_i H_i}=\{ 0\}$;
\item[(6)] $\DD_{\{0\},H}=\{ 0\}$ for all $H$;
\item[(7)] $\DD_{E,\R} =\{ 0\}$ for all $E$;
\item[(8)] $\DD_{\Z t,H}=\{ 0\}$ for all $t\in \R$ and all $H$;
\item[(9)] if $E+H$ is admissible, then $\DD_{E,F}=\{ 0\}$; in particular, if $E$ and $H$ are countable, then
$\DD_{E,H}=\{ 0\}$.
\end{itemize}

Property (6) is a consequence of the fact that $E=\{0\}$ is representable
(it is represented for example by $B=\PP$), (7) is obvious, (8) holds since the cyclic subgroups of $\R$
are representable. The first part of (9) is a mere consequence of the definition of $\DD_{E,H}$. The second part of (9) follows
from Corollary 5.16.

\begin{lem5.21*}
\emph{Let $E$ and $H$ be subgroups of $\R$ such that $E\subset G_0$ for some representable $G_0$. Then $H^\prime =\bigcap_{G\in\RR_E} (G+H)$
is the smallest subgroup of $\R$ which contains $H$ and satisfies $\DD_{E,H^\prime} =\{ 0\}$.
Moreover, $\DD_{E,H}=H^\prime/H$.}
\end{lem5.21*}

\begin{proof}
One plainly checks that $\bigcap_{G\in\RR_E} (G+H^\prime)=H^\prime$. Since $H^\prime$ contains $E$, we have
$E+H^\prime=H^\prime$ and $\DD_{E,H^\prime}=\{ 0\}$. If $H\subset H_1$ and $\DD_{E,H_1}=\{ 0\}$, then
\begin{equation*}
H^\prime =\bigcap\limits_{G\in\RR_E} (G+H) \subset \bigcap\limits_{G\in\RR_E} (G+H_1) .\qedhere
\end{equation*}
\end{proof}

\begin{lem5.22*}
\emph{If $E$ and $H$ are subgroups of $\R$ such that $E$ is contained in a representable group $G_0$, then
$\DD_{E,H} \simeq \DD_{E,H\cap G_0}$, the isomorphism being induced by the inclusion $H\cap G_0 \subset H$.}
\end{lem5.22*}

\begin{proof}
Denote
\begin{equation*}
\begin{split}
G_1 & = \{ G+(H\cap H_0): G\in \RR_E, G\subset G_0 \}, \\
G_2 & =\{ G+H: G\in \RR_E, G\subset G_0 \} .
\end{split}
\end{equation*}

We have $H^\prime =\bigcap_{G\in\RR_E} (G+H)=G_2$ and the similar equality with $H\cap G_0$ instead of $H$ and $G_1$ instead of $G_2$.
One implication is obvious. The other one follows from the fact that the intersection of two representable groups is representable.

It is simple to check that the natural morphism $G_1 \rightarrow G_2/(E+H)$ induced by the identity is surjective,
for $G_2=E+H+G_1$ and its kernel is $E+(H\cap G_0)$. We therefore obtain
\begin{equation*}
\DD_{E,H} =G_2 /(E+H) \simeq G_1 /(E+(H\cap G_0)) =\DD_{E,H\cap G_0} .\qedhere
\end{equation*}
\end{proof}

\begin{def5.23*}
Let $E$, $H_1$ and $H_2$ be subgroups of $\R$ such that $E$ is contained in a representable group. We say that
$H_1$ \emph{is equivalent to $H_2$ modulo $E$} (write $H_1 \sim_E H_2$) if there exists a representable group $G$ containing
$E$ such that $G\cap H_1=G \cap H_2$.
\end{def5.23*}

\begin{cor5.24*}
\emph{If $H_1 \sim_E H_2$, then $\DD_{E,H_1} \simeq \DD_{E,H_2}$.}
\end{cor5.24*}

\begin{lem5.25*}
\emph{Let $E$ be a countable subgroup of $\R$ and $G_n$ be a countable family of representable groups such that
$E\subset G =\bigcap_n G_n$. Then, there exists a representable group $\Gamma$ such that $E \subset \Gamma \subset G$.}
\end{lem5.25*}

\begin{proof}
Let $E=\bigcup_n E_n$, with $\{ E_n\}_n$ an increasing sequence of finite subsets of $E$. Assume that $G_n=G(f_{B_n})$, with
$B_n$ subsets of $\PP$. Since $f_{B_n}(t)<\infty$ for all $t\in E$ and all $n$, we may choose for each $n$ a finite subset $F_n$ of $B_n$ such that
if we take $B_n^\prime =B_n \setminus F_n$, then
\begin{equation*}
f_{B_n^\prime}(t) < 2^{-n},\quad \mbox{\rm for $t\in F_n$.}
\end{equation*}

We set $B^\prime =\bigcup_n B_n^\prime$. Then we get
\begin{equation*}
f_{B^\prime} (t)=\sum\limits_{j\geqslant 1} f_{B_j^\prime}(t) \leqslant\sum\limits_{1\leqslant j<k} f_{B_j^\prime} (t) +\sum\limits_{j\geqslant k} 2^{-j} < \infty ,
\end{equation*}
and consequently $E\subset G(f_{B^\prime}) \subset \bigcap_k G(f_{B_k^\prime}) =\bigcap_k G(f_{B_k}) =G$.
\end{proof}

\begin{cor5.26*}
\emph{Let $E$ and $H$ be subgroups of $\R$ such that $E$ is countable. Let $\{ G_n\}_n$ be a countable family of representable groups
such that $E\subset \bigcap_n G_n$ and define $\{ H_n\}_n$ by $H_1 =H$, $H_{n+1}=H_n \cap G_n$ if $n\geqslant 1$.
Then, the natural group morphism $\DD_{E, \cap_n H_n} \rightarrow \DD_{E,H}$ is an isomorphism.}
\end{cor5.26*}

\medskip

The main results of this paper were announced in the note `Facteurs de type III associ\' es aux ensembles de nombres premiers',
\emph{C. R. Acad. Sci. Paris S\' er. I} 324 (1997), 797--800.

\bigskip

\emph{Acknowledgments.}
We are grateful to U. Haagerup, G. Skandalis and M. Takesaki for bringing to our attention references
\cite{CW}, \cite{GS} and \cite{Su}.

\end{document}